\documentclass[journal]{new-aiaa}
\usepackage[utf8]{inputenc}
\usepackage{textcomp}

\usepackage{graphicx}
\usepackage{amsmath}
\usepackage[version=4]{mhchem}
\usepackage{siunitx}
\usepackage{longtable,tabularx}
\usepackage{booktabs}
\title{Joint Replenishment Strategy for Multiple Satellite Constellations with Shared Launch Opportunities}

\author{Jaewoo Kim\footnote{Ph.D. Candidate, Department of Aerospace Engineering, 291 Daehak-Ro.}}
\affil{Korea Advanced Institute of Science and Technology, Daejeon 34141, Republic of Korea}
\author{Taehyun Sung\footnote{Senior Researcher}}
\affil{Agency for Defense Development, Daejeon 34186, Republic of Korea}
\author{Woonam Hwang\footnote{Associate Professor of Industrial and Systems Engineering, 291 Daehak-Ro.}}
\affil{Korea Advanced Institute of Science and Technology, Daejeon 34141, Republic of Korea}
\author{Jaemyung Ahn\footnote{Professor of Aerospace Engineering, 291 Daehak Ro; jaemyung.ahn@kaist.ac.kr. Associate Fellow AIAA (Corresponding Author).}}
\affil{Korea Advanced Institute of Science and Technology, Daejeon 34141, Republic of Korea}
\begin{document}

\maketitle

\begin{abstract}
This paper proposes a novel replenishment strategy that can jointly support multiple satellite constellations. In this approach, multiple constellations share launch opportunities and parking orbits to address the operational satellite failures and ensure the desired service level of the constellations. We develop an inventory management model based on parametric replenishment policies, considering the launch vehicle's capacity and the shipping size of satellites. Based on this model, we introduce two decision-making scenarios and propose their corresponding solution frameworks. We conduct two case studies to provide valuable insights into the proposed strategy and demonstrate its applicability to supply chain management for maintaining multiple satellite constellations.
\end{abstract}

\section*{Nomenclature}
{\renewcommand\arraystretch{1.0}
\noindent\begin{longtable*}{@{}l @{\quad=\quad} l@{}}
$A$ & limited shipping capacity of a vehicle, in units of shipping sizes\\
$\mathbf{a}$ & batch volume vectors, in units of shipping sizes\\
$C_{\text{hold}}$ & annual holding cost, \$M\\
$C_{\text{lau}}$ & annual launch cost, \$M\\
$C_{\text{maneuv}}$ & annual maneuvering cost, \$M\\
$C_{\text{manufac}}$ & annual manufacturing cost, \$M\\
$c_{\text{hold}}$ & annual holding cost of a satellite, \$M/satellite/year\\
$c_{\text{lau}}$ & launch costs, \$M\\
$c_{\text{manufac}}$ & manufacturing cost of a satellite, \$M/satellite\\
$EQ_{\text{parking}}$ & expected order size of batches at parking spare, per unit time\\
$ES_{\text{parking}}$ & expected shortage of parking spares, -\\
$ES_{\text{plane}}$ & expected shortage of in-plane spares, in units of satellites\\
$\mathbf{e}$ & vector of which all the components are ones\\
$\mathbf{e}_j$ & vector of which $j^{\text{th}}$ component is one but zeros\\
$h_{\text{parking}}$ & altitude of parking orbits, km\\
$h_{\text{plane}}$ & altitude of orbital planes in constellation, km\\
$i$ & inclination angle of an orbit, deg \\
$M_{\text{dry}}$ & dry mass of satellite, kg\\
$M_{\text{fuel}}$ & fuel mass required for orbital transfer, kg\\
$\dot{M}_{\text{prop}}$ & mass flow rate of satellite propulsion system, kg/s\\
$m$ & number of satellite constellations to jointly replenish\\
$N_{\text{parking}}$ & number of parking orbits, in units of planes\\
$N_{\text{plane}}$ & number of orbital planes in constellation, in units of planes\\
$N_{\text{sat}}$ & number of satellites per orbital plane in constellation, in units of satellites\\
$N_t$ & number of time units per year\\
$n_{\text{parking}}$ & order frequency of parking spares, per unit time\\ 
$n_{\text{plane}}$ & order frequency of in-plane spares, per unit time\\
$Q$ & batch size of in-plane spares, in units of satellites \\
$\overline{SL_{\text{parking}}}$ & mean stock level of parking spares, -\\
$\overline{SL_{\text{plane}}}$ & mean stock level of in-plane spares, in units of satellites\\ 
$s$ & reorder point of in-plane spares, in units of satellites \\
$S$ & order-up-to level of parking spares, -\\
$t_{\text{trans}}$ & in-plane transfer maneuver time, in units of time\\
$T_{\text{prop}}$ & thrust of satellite propulsion system, N \\
$\mathcal{T}_{n}$ & admissible time interval for supplementation from the n$^{\text{th}}$ closest parking orbit\\
$U$ & size reorder point of parking spare, in units of shipping sizes \\
$V_{\text{ex}}$ & exhaust velocity of satellite propulsion system, km/s\\
$\mathbf{w}$ & inventory state vector\\
$\beta$ & proxy for bargaining power\\ 
$\delta$ & stochastic demand \\
$\epsilon_{\text{maneuv}}$ & fuel mass conversion coefficient, \$M/kg\\ 
$\lambda_{\text{parking}}$ & demand rate at parking spares, per unit time \\
$\lambda_{\text{plane}}$ & demand rate at in-plane spares, in units of satellites per unit time \\
$\lambda_{\text{sat}}$ & failure rate of a satellite, failures per year \\
$\rho_{\text{parking}}$ & order fill rate for parking spares, -\\
$\rho_{\text{plane}}$ & order fill rate for in-plane spares, -\\ 
$\tau$ & stochastic lead time\\
$\omega$ & stochastic inventory state\\
$\mathbb{R}_+$ & non-negative real number set \\
$\mathbb{Z}_+$, $\mathbb{Z}_{++}$ & non-negative, and positive integer sets\\
$\lfloor\cdot\rfloor$ & floor operator returning the greatest integer less than or equal to the input variable\\
OUL & order-up-to level\\
ROP & reorder point\\
SROP & shipping size reorder point\\
\multicolumn{2}{@{}l}{Subscripts}\\
$j$ & index for constellation\\
$k$ & index for parking orbit\\
\end{longtable*}}

\section{Introduction}\label{sec: section 1}
\lettrine{T}{he} space industry is experiencing unparalleled opportunities, driven by increasing demands for global coverage services and technological advancements. Innovations such as the miniaturization of electronics, cost-effective satellite manufacturing using Commercial Off-The-Shelf (COTS) components, and reduced launch service costs facilitated by reusable launch vehicles and shared launch opportunities enable the widespread deployment of small, low-cost satellites \cite{poghosyan2016unified}. Satellite constellations are a key part of this trend, and these systems are now being scaled extensively in Low Earth Orbit (LEO) to address diverse market needs. Various mega-constellation initiatives, such as Eutelsat's OneWeb, SpaceX's Starlink, and Amazon's Project Kuiper, exemplify this trend \cite{pachler2021updated, pachler2024flooding}.

Satellite mega-constellations pose unique challenges for system maintenance compared to traditional satellite operations. Their vast scale and reliance on less reliable, low-cost satellites require operators to address more frequent satellite failures. This situation necessitates a rethinking of maintenance strategies, including decisions on where to position spare satellites, when and how much to replenish spare inventory at each level, and how to interconnect and transfer spares between different levels---the \textit{logistics} of spare satellites.

Strategic approaches for maintaining mega-constellations are just beginning to emerge in the literature. Jakob et al. \cite{jakob2019optimal} developed a multi-echelon inventory management model specifically for maintaining a satellite mega-constellation. This model features a three-echelon inventory system comprising ground spares, parking spares, and in-orbit spares, which have terrestrial analogies to suppliers, warehouses, and retailers, respectively. Notably, parking orbits located at altitudes different from the operational orbital planes act as shared warehouses for in-plane spares. Within this model, the parking spares and in-plane spares are replenished by different $\left(s,Q\right)$ policies \cite{silver2016inventory}, where a batch of $Q$ units is ordered when the stock level drops to $s$. They formulated analytical expressions for inventory-related parameters, maintenance costs, and service level, assuming the Poisson distribution for the failure and the exponential distribution for the launch lead time. The costs encompass satellite manufacturing, launch, holding, and maneuvering from parking orbits to operational orbits. They introduced an optimization problem for identifying the optimal policy parameters and parking orbit geometry that minimize the maintenance costs while ensuring the service level for spare inventory.

Kim et al. \cite{kim2025replenishment} extended this multi-echelon inventory model and proposed a strategy with dual supply modes, incorporating a direct channel that injects the spare satellites into the orbital planes \cite{sung2023optimal}. They integrated an inventory policy denoted as $\left(R_1,R_2,Q_1,Q_2\right)$ with a time window for the in-plane spares. In this policy, when the stock level drops to $R_1$, $Q_1$ units are ordered from the parking spares, the primary supply mode (indirect channel). If the stock level drops to $R_2$ within a predefined time window, the auxiliary supply mode is activated (direct channel), and $Q_2$ units are ordered. They developed analytic expressions for inventory-related parameters and introduced two optimization problems from different perspectives---the satellite constellation operator and the auxiliary launch service provider.

Although research on strategic initiatives remains limited, the space industry's rapid evolution continues to reveal new opportunities. For example, as multiple satellite constellation projects emerge, operators can strategically partner to reduce costs---an approach widely used in terrestrial industries. This is particularly compelling as launch vehicle capacities increase, making shared launch missions more efficient than exclusive use. Yet this collaborative scenario remains underexplored in current literature.

To address the identified research gaps, we propose a novel framework for the joint replenishment strategy of multiple satellite constellations. Specifically, we develop an inventory management model to support decision-making during the strategic or planning phase of a satellite mega-constellation project. In this model, one authority or multiple authorities replenish multiple satellite constellations by shared launch opportunities and parking orbits. We introduce relevant decision-making scenarios and corresponding solution procedures based on this model.

The rest of this paper is organized into five sections. Section \ref{sec: section 2} introduces critical concepts related to satellite constellations and the fundamental principles for developing the mathematical model. Section \ref{sec: section 3} describes the supply chain for the joint replenishment strategy of multiple satellite constellations, constructs the corresponding mathematical model, and validates the model. Section \ref{sec: section 4} explains the decision-making contexts associated with the joint replenishment strategy and the solution procedures. Section \ref{sec: section 5} presents case studies to illustrate scenarios where the joint replenishment strategy proves efficient and demonstrates the applicability of the proposed approach. Finally, Section \ref{sec: section 6} concludes the paper.

\section{Preliminaries}\label{sec: section 2}
This section outlines the fundamental concepts and establishes the language. It first introduces satellite constellations and orbital mechanics, then presents relevant concepts from inventory management.

\subsection{Satellite Constellation and Orbital Mechanics}
This subsection describes the Walker pattern of satellite constellation, $J_2$ perturbation, and orbital transfer \cite{wertz2001mission, edberg2020design, vallado2022fundamentals}.

\subsubsection{Walker Delta Pattern}
A constellation is ``a set of satellites distributed over space to work together to achieve common objectives'' \cite{wertz2001mission}. Various design methodologies have been proposed, and the \textit{Walker Delta pattern} \cite{walker1970circular} is one of them. A Walker Delta pattern is characterized by four parameters, $i:t/p/f$. Here, $t$ denotes the number of satellites evenly distributed among $p$ orbital planes. Each orbital plane contains $s$ satellites, and all the orbital planes share the common inclination angle $i$. The right ascension of the ascending node (RAANs) of each orbital plane is spaced along the equator, with uniform intervals of $2\pi/p$. Satellites are located at intervals of $2\pi/s$ within each plane. The phase offset between satellites in neighboring orbital planes is determined as an integer multiple $f$ of $2\pi/t$. This symmetrical arrangement ensures that all satellites are subjected to similar orbital perturbations. This work assumes that the satellite constellations adhere to Walker Delta patterns, with $N_{\text{plane}}$ and $N_{\text{sat}}$ substituting $p$ and $s$, respectively. Additionally, we assume all the orbital planes are circular with altitudes of $h_{\text{plane}}$.

\subsubsection{\texorpdfstring{$J_2$ Perturbation and RAAN drift}{}}
The $J_2$ effect is considered the sole perturbing factor influencing satellite motion in this context. This perturbation induces long-term changes in three orbital elements: the RAAN, the argument of perigee, and the mean anomaly. Since circular orbits are assumed in this study, the primary perturbation of interest is RAAN drift. This drift causes the RAAN of an orbital plane to gradually shift westward or eastward over time at a specific rate ($\dot{\Omega}$). The rate of this nodal precession is determined depending on several factors as follows \cite{vallado2022fundamentals}:
\begin{equation}\label{eqn: nodal precession}
    \dot{\Omega}\left(a,e,i\right) =-\frac{3}{2} \sqrt{\frac{\mu_E}{a^3}} \frac{R_E^2}{a^2 (1-e^2)^2} J_2 \cos{i},
\end{equation}
where $a$, $i$, and $e$ are the semi-major axis, the inclination, and the eccentricity of the orbit, $J_2$ is the zonal harmonic coefficient, $R_E$ is the Earth's equatorial radius, and $\mu_E$ represents the Earth's gravitational constant.

\subsubsection{Orbital Transfer}
A satellite transfers from one orbit to another using its propulsion system. This work considers only transfers between coplanar circular orbits under an inverse-square gravitational field with finite continuous thrust. Low-thrust propulsion systems, such as electric propulsion, have been widely applied for continuous-thrust maneuvers of small satellites due to their high efficiency. Alfano et al. \cite{alfano1994constant} modeled the velocity increment required for continuous low-thrust transfers between circular orbits, and Bettinger et al. \cite{bettinger2014mathematical} demonstrated the equivalence of the velocity increments for Hohmann transfers and continuous low-thrust transfers. Following these works, the required velocity increment ($\Delta V$) and the time of flight ($TOF$) for continuous low-thrust maneuvers are given as
\begin{equation}\label{eqn: delta v}
    \Delta V = \sqrt{\frac{\mu_E}{r_i}} - \sqrt{\frac{\mu_E}{r_f}},
\end{equation}
\begin{equation}\label{eqn: time of flight}
    TOF = \frac{M_{\text{fuel}}}{\dot{M}_{\text{prop}}},
\end{equation}
where $r_i$ and $r_f$ represent the initial and final orbit radii ($r_i<r_f$), respectively, $\mu_E$ is Earth's gravitational constant, $M_{\text{fuel}}$ denotes the required fuel mass, and $\dot{M}_{\text{prop}}$ is the mass flow rate of the propulsion system.

The velocity increment ($\Delta V$) to the required fuel mass ($M_{\text{fuel}}$) are related by the rocket equation 
\begin{equation}\label{eqn: rocket equation}
    M_{\text{fuel}} = M_{\text{dry}}(e^{\Delta V/V_{\text{ex}}}-1),
\end{equation}
where $M_{\text{dry}}$ and $V_{\text{ex}}$ are the satellite's dry mass and the exhaust velocity of the propulsion system, respectively.

\subsection{Inventory Management}
Inventory management systematically controls inventory to optimize supply chain performance. Effective inventory management primarily focuses on cost-effectiveness and fulfilling customer demand, and the former is the emphasis of this research. The operator utilizes several levers to achieve this, including determining when to place orders, deciding the order quantity, and selecting the appropriate transportation. This section outlines the core inventory management policies with stochastic lead time and demand, which are relevant to the timing and quantity of an order. It begins with an introduction to fundamental terminology, followed by a discussion of the $\left(s, Q\right)$ and $\left(U, \mathbf{S}\right)$ policies.

\subsubsection{Representation of Inventory Status}
We use the following concepts for representing the inventory in the proposed model \cite{silver2016inventory}:
\begin{itemize}
    \item On-hand stock ($H$) is the stock that is physically on the shelf; in our context, this refers to spare satellites ready for resupply to the lower echelon. Consequently, this value is always non-negative.
    \item Backorders ($B$), also called shortages, represent the demands not fulfilled by on-hand stock. These unmet demands are satisfied by the next order. This value is always non-negative.
    \item Net stock ($NS$) is defined by the following relation:
    \begin{equation}\label{eqn: net stock definition}
        NS = H - B.
    \end{equation}
    This can have a negative value. Since backorders do not occur when on-hand stock can fulfill all the demands, the following relations hold:
    \begin{equation}\label{eqn: net stock additional relations}
        \text{max}\left\{NS,0\right\} = H,\ \text{max}\left\{-NS,0\right\} = B.
    \end{equation}
    \item On-order stock ($O$) is the stock that has been ordered but has not yet been received.
    \item Inventory position ($IP$) is the sum of the net stock and the on-order stock, represented by the following equation:
    \begin{equation}\label{eqn: inventory position definition}
        IP = NS + O = H - B + O.
    \end{equation}
    This work often refers to the inventory position as the stock level.
\end{itemize}

\subsubsection{\texorpdfstring{$\left(s,Q\right)$}{} Policy}\label{section: (s, Q) policy}
The $\left(s,Q\right)$ policy is a continuous review system for inventory management, meaning the stock level is reviewed and replenished in a continuous manner \cite{silver2016inventory}. Under this policy, a fixed order quantity ($Q$) is ordered whenever the inventory position reaches the reorder point (ROP), $s$. Figure \ref{fig: (s, Q) policy in order cycle representation} illustrates an order cycle following an $\left(s,Q\right)$ policy. With this policy, inventory parameters, such as the mean stock level, expected shortage within an order cycle, and order frequency, can be obtained as closed forms.
\begin{figure}[hbt!]
\centering
\includegraphics[page=1, width=0.7\textwidth]{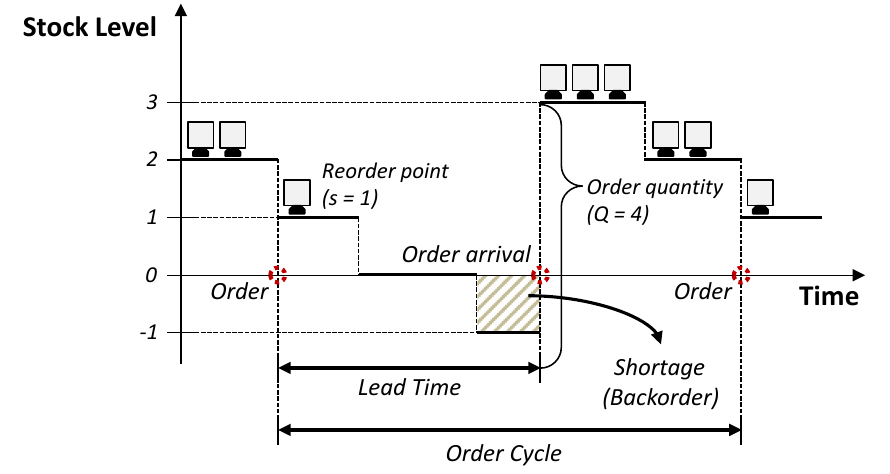}
\caption{Order cycle representation with an $\left(s,Q\right)$ policy ($s=1$, $Q=4$)}
\label{fig: (s, Q) policy in order cycle representation}
\end{figure}

The mean stock level ($\overline{SL}$) represents the average stock maintained at the operator's facility. Under the assumption that shortages are rare, the net stock can replace on-hand stock. Consequently, the expected on-hand stock immediately before an order arrives ($\mathbb{E}[H_1]$) is defined as
\begin{equation}
    \mathbb{E}[H_1] = s - \mathbb{E}\left[x_\tau\right],
\end{equation}
where $x_\tau$ is the demand during the lead time ($\tau$). Since the order quantity is $Q$, the expected on-hand stock just after the arrival of an order ($\mathbb{E}[H_2]$) is given as
\begin{equation}
    \mathbb{E}[H_2] = \mathbb{E}[H_1] + Q = s - \mathbb{E}[x_\tau] + Q.
\end{equation}
With a linear approximation in the stock level drop, the mean stock level ($\overline{SL}$) is
\begin{equation}
    \overline{SL} = \frac{\mathbb{E}[H_1] + \mathbb{E}[H_2]}{2} + \frac{1}{2} = s - \mathbb{E}[x_\tau] + \frac{Q}{2} + \frac{1}{2},
\end{equation}
where the term $1/2$ is the continuity correction factor, which adjusts the discrepancy between the assumption of linear stock drops and the reality of discrete drops. When demand follows a Poisson process with rate $\lambda$, the mean stock level during lead time $\tau$ is represented as
\begin{equation}\label{eqn: mean stock level Poisson}
    \overline{SL}\left(s, Q; \lambda, \tau\right) = s - \lambda \tau + \frac{Q}{2} + \frac{1}{2}.
\end{equation}

With an $\left(s, Q\right)$ policy, the expected shortage ($ES$) is the mean of the excessive demand over the ROP ($s$) within an order cycle. If demand during the lead time follows a probability mass function (PMF) $P_{D}$, the general expression for expected shortage per order cycle is 
\begin{equation}\label{eqn: expected shortage integral form}
    ES\left(s\right) = \sum_{\delta=s}^{\infty} \left( \delta - s \right) P_{D}\left(\delta\right).
\end{equation}
In the specific case where demand follows a Poisson process with a demand rate of $\lambda$, the expected shortage during the lead time $\tau$ is
\begin{equation}\label{eqn: expected shortage summation form with Poisson}
    ES\left(s; \lambda, \tau\right) = \sum_{\delta=s}^{\infty} \left( \delta - s \right) \frac{ e^{-\lambda \tau} \left( \lambda \tau \right)^\delta }{\delta!}.
\end{equation}
Expanding Eq. \eqref{eqn: expected shortage summation form with Poisson} leads to
\begin{equation}\label{eqn: expected shortage expansion with Poisson}
    ES\left( s; \lambda, \tau \right) = \lambda \tau \left( 1 - F_{\lambda}(s - 1;\tau) \right) - s \left( 1 - F_{\lambda}(s;\tau) \right),
\end{equation}
where $F_{\lambda}\left(\,\cdot\,;\tau\right)$ is the cumulative distribution function (CDF) of the Poisson distribution with parameter $\lambda \tau$ \cite{ganeshan1999managing}. Since satellite failures are modeled as a Poisson process in this work, Eq. \eqref{eqn: expected shortage expansion with Poisson} is used for calculating the expected shortage. The order fill rate ($\rho$) measures the average fraction of demand immediately fulfilled by the on-hand stock, which is given as
\begin{equation}\label{eqn: order fill rate}
    \rho = 1 - \frac{ES}{Q},
\end{equation}
where $ES$ is the expected shortage, and $Q$ is the fixed order quantity.

In addition, the order frequency ($n$), which indicates the expected number of orders (or order cycles) per unit time, is defined as
\begin{equation} \label{eqn: order frequency}
    n = \frac{\lambda}{Q}.
\end{equation}

\subsubsection{\texorpdfstring{$\left(U, \mathbf{S}\right)$}{} Policy}\label{section: (U, S) policy}
The $\left(U, \mathbf{S}\right)$ policy is a parametric replenishment strategy for the joint replenishment problem (JRP), where $m$ different stock-keeping units (SKUs) are replenished together through joint ordering. Since different SKUs have distinct characteristics (e.g., varying demand rates), it is crucial for inventory management to address these differences effectively.

Although some past studies focused on the JRP \cite{khouja2008review, bastos2017systematic}, research that considers different shipping sizes and stochastic demands---factors essential for our context of maintaining satellite mega-constellations---remains limited. Only a few studies have explored these issues within the JRP framework \cite{cachon2001managing, mustafa2010joint, li2020stochastic}, with notable contributions from Li and Schmidt \cite{li2020stochastic}. Their work studied replenishing SKUs with varying shipping sizes through joint ordering under a limited shipping capacity ($A$). Each SKU has a shipping size of $a_j$ units with a demand following a Poisson process with rate parameter $\lambda_j$ ($j = 1, 2, \dots, m$). By setting the order-up-to level (OUL) $S_j$ for each SKU, inventory is replenished up to these OULs through a joint order. This joint order is triggered whenever the total shipping size of cumulative demand awaiting replenishment reaches $U$ ($\leq A$), the shipping size reorder point (SROP).

Li and Schmidt \cite{li2020stochastic} adopted the Markov chain (MC) as their foundational analytical framework, where $\mathbf{w}$ ($\in\mathbb{Z}_{+}^m$) represents the inventory state, defined as
\begin{equation} \label{eqn: state definition}
    \mathbf{w}=\left[w_1,w_2,\dots,w_m\right]^\top=\left[S_1-IP_1,S_2-IP_2,\dots,S_m-IP_m\right]^\top,
\end{equation}
where $IP_j$ is the inventory position of SKU $j$. Accordingly, each component of $\mathbf{w}$ represents the stock level drop from the OUL of the corresponding SKU, and a demand arrival of SKU $j$ increments the $j$th component of the state vector ($\mathbf{w}$).

Note that under the $\left(U,\mathbf{S}\right)$ policy, any state $\mathbf{w}$ satisfying $\mathbf{a}^\top\mathbf{w}\ge U$ cannot be reached because a joint order is triggered when the cumulative shipping size of demand ($\mathbf{a}^\top\mathbf{w}$) touches the SROP ($U$). Therefore, the state space is defined as:
\begin{equation} \label{eqn: state space}
    \mathcal{S}_0=\left\{\mathbf{w}\in \mathbb{Z}_+^m: \mathbf{a}^\top \mathbf{w}< U\right\}.
\end{equation}

Given this state space, a demand arrival causes a state transition. There are three potential outcomes when demand for SKU $j$ arrives:
\begin{enumerate}
    \item Outcome 1: If $\mathbf{w}+\mathbf{e}_j\in \mathcal{S}_0$, an order is not triggered, and the state progresses to $\mathbf{w}'=\mathbf{w}+\mathbf{e}_j$.
    \item Outcome 2: If $\mathbf{w}+\mathbf{e}_j\notin \mathcal{S}_0$ and $\mathbf{w}+\mathbf{e}_j\le A$, an order is triggered, and all cumulative demands are replenished. The state proceeds to $\mathbf{w}'=\mathbf{0}$.
    \item Outcome 3: If $\mathbf{w}+\mathbf{e}_j\notin \mathcal{S}_0$ and $\mathbf{w}+\mathbf{e}_j> A$, an order is triggered, but not all the demands can be covered by the order, leaving the last demand ($\mathbf{e}_j$) unmet and waiting for the next order. Consequently, the state transitions to $\mathbf{w}'=\mathbf{e}_j$.
\end{enumerate}
Here, $\mathcal{I}_\mathbf{w}$ denotes the index set of SKUs that trigger Outcome 2 given the current inventory state $\mathbf{w}$, defined as:
\begin{equation}\label{eqn: I_w definition}
    \mathcal{I}_{\mathbf{w}}=\left\{j\in \left\{1,2,\dots,m\right\}:\mathbf{w}+\mathbf{e}_j\notin \mathcal{S}_0,\ \mathbf{a}^\top \left(\mathbf{w}+\mathbf{e}_j\right)\le A\right\}.
\end{equation}

\begin{figure}[hbt!]
\centering
\includegraphics[page=2, width=0.7\textwidth]{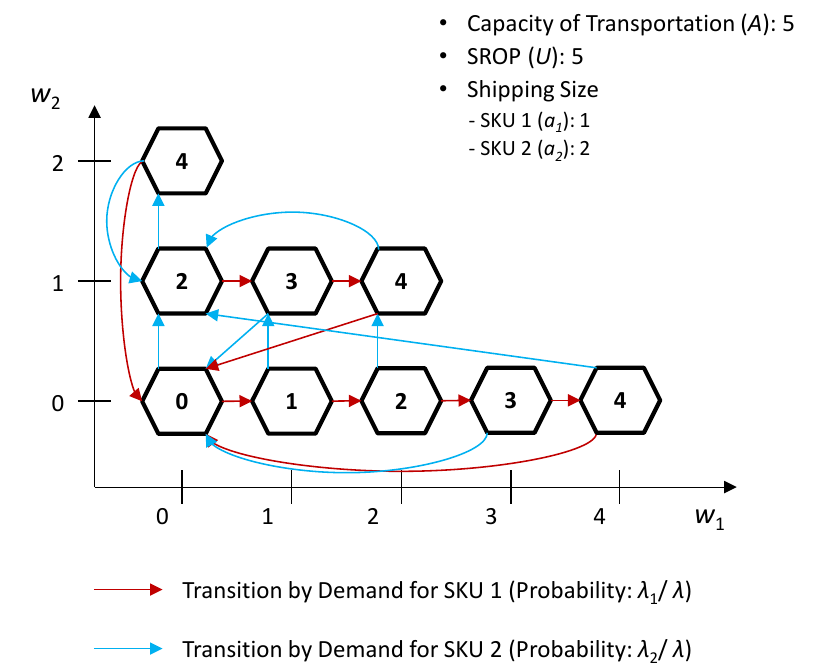}
\caption{Example MC: The number of SKUs is 2 ($m = 2$) with shipping sizes of 1 ($a_1=1$) and 2 ($a_2=2$) for each SKU; both the SROP ($U$) and the transportation capacity ($A$) are set to 5}
\label{fig: markov chain example}
\end{figure}

Since the demand for each SKU $j$ follows a Poisson process with parameter $\lambda_j$ and is independent of the others, the probability that $\mathbf{w}$ changes due to the arrival of demand for SKU $j$ is $\lambda_j / \lambda$, where $\lambda = \sum_{j=1}^m \lambda_j$.
The transition probability of the MC representing the order cycle is defined as
\begin{equation}\label{eqn: state transition probability}
    p_{\mathbf{w}\mathbf{w}'} = \begin{cases}
        \lambda_{j}/\lambda,& \text{when }\mathbf{w}'=\mathbf{w}+\mathbf{e}_j\text{ or } \mathbf{w}'=\mathbf{e}_j\text{ and }\mathbf{a}^\top \left(\mathbf{w}+\mathbf{e}_j\right)>A \text{ (Outcome 1 \& 3)},\\
        \sum_{j\in \mathcal{I}_{\mathbf{w}}}\lambda_j/\lambda,& \text{when }\mathbf{w}'=\mathbf{0} \text{ (Outcome 2)},\\
        0,&\text{otherwise}.
    \end{cases}
\end{equation}
An example MC with two SKUs is illustrated in Fig. \ref{fig: markov chain example}. In this figure, each hexagon contains a number representing the shipping size of the cumulative demand from the OULs.

If an MC is irreducible and has a finite state space, the stationary distribution vector $\mathbf{w}_\infty$ ($\in \mathbb{R}_{+}^{\left|\mathcal{S}\right|}$) represents the stationary distribution when the chain is aperiodic, or at least the limiting frequency in cases where the MC exhibits periodicity \cite{rosenthal2006first}. The stationary distribution vector $\mathbf{w}_\infty$ is defined as the vector satisfying 
\begin{equation}\label{eqn: stationary distribution condition 1} P^\top\mathbf{w}_{\infty} = \mathbf{w}_{\infty}, \quad \text{where } P=\left(p_{\mathbf{w},\mathbf{w}'}\right), \quad \mathbf{w}, \mathbf{w}'\in \mathcal{S}_0;\ \text{and}
\end{equation}
\begin{equation}\label{eqn: stationary distribution condition 2} \sum_{\mathbf{w}\in\mathcal{S}_0}{w_{\infty,\mathbf{w}}}=1,
\end{equation}
where $w_{\infty,\mathbf{w}}$ denotes the component of $\mathbf{w}_\infty$ corresponding to the state $\mathbf{w}$. This stationary distribution describes the mean behavior of the inventory in the infinite time horizon.

Li and Schmidt \cite{li2020stochastic} provided detailed derivations of various inventory parameters. This work extends their model by incorporating stochastic lead time and applying it to model the inventory of parking spares.
 
\section{Model Construction}\label{sec: section 3}
This section outlines the proposed supply chain and the corresponding inventory management model. It begins with an overview of the supply chain and the assumptions used to construct the model. Next, it presents a detailed mathematical model for the inventory behavior at each echelon, connecting it to the total maintenance costs. Lastly, we evaluate the model's accuracy by comparing its results with those obtained from Monte Carlo simulations.

\subsection{Model Overview and Assumptions}
We propose a joint replenishment strategy for $m$ different constellations. Each constellation has a distinct Walker Delta pattern. Constellation $j$ ($=1,2,\dots,m$) have $N_{\text{plane},j}$ circular orbital planes with an altitude of $h_{\text{plane},j}$ and $N_{\text{sat},j}$ operational satellites. The RAANs of the orbital planes of the constellation $j$ are distributed at a regular interval of $2\pi/N_{\text{plane},j}$, and all the constellations share a common inclination angle $i$. As a result, the model permits variations among the constellations only in altitude ($h_{\text{plane},j}$), the number of orbital planes ($N_{\text{plane},j}$), and the number of satellites per plane ($N_{\text{sat},j}$).

Different constellations may consist of distinct satellites characterized by their dry mass ($M_{\text{dry},j}$), shipping size ($v_{\text{j}}$), and failure rates ($\lambda_{\text{sat},j}$). Here, ``shipping size'' refers to the number of docking plate or adapter slots that a single satellite occupies \cite{spacex2024rideshare}. We assume that the failure rate of each satellite is time-invariant and model failures using a Poisson process.

\begin{figure}[hbt!]
    \centering
    \includegraphics[page=3, width=1.0\textwidth]{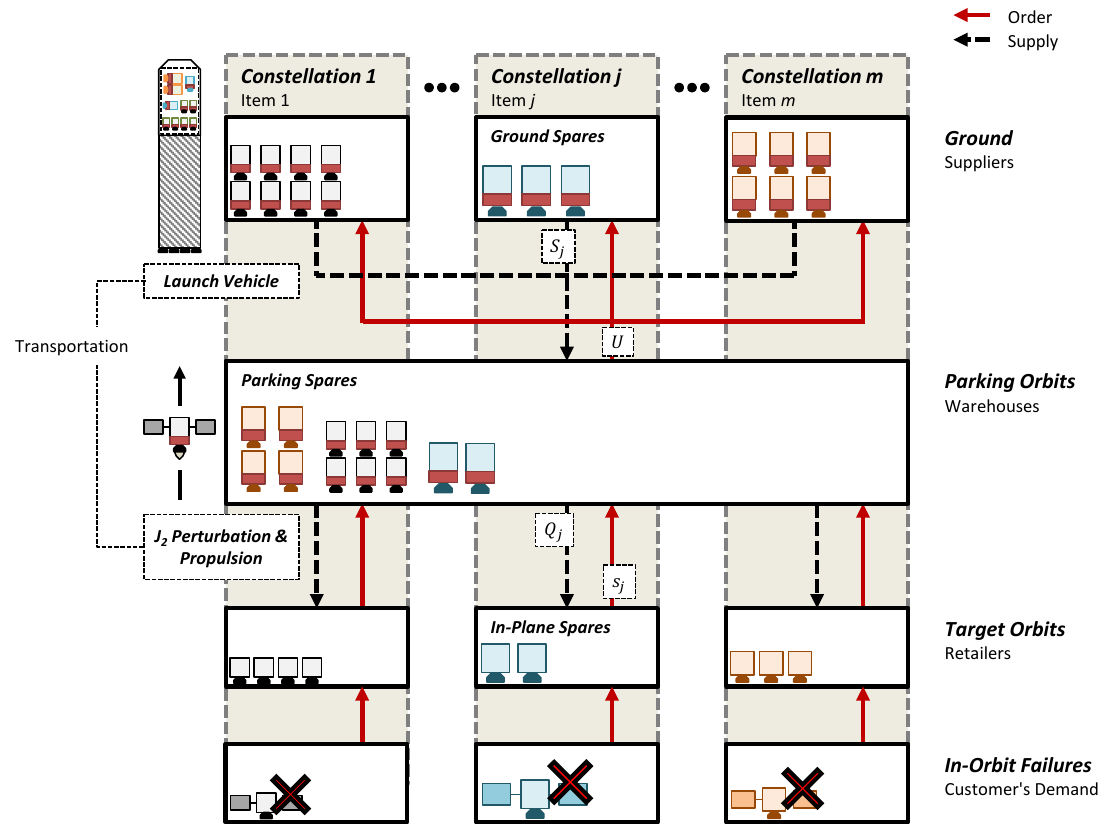}
    \caption{Overview of the proposed inventory model with ground analogy}
    \label{fig: overview of the supply chain for joint replenishment strategy}
\end{figure}

The proposed supply chain consists of three echelons of spare satellites: \textit{in-plane spares}, \textit{parking spares}, and \textit{ground spares}. Figure \ref{fig: overview of the supply chain for joint replenishment strategy} illustrates the proposed supply chain with a ground analogy.

The first echelon, the \textit{in-plane spares}, immediately responds to the in-orbit failures of the satellites. Each orbital plane has its in-orbit spares and sources the spare satellites from them. We assume the replacement time using the in-orbit spares (1--2 days) is negligible \cite{jakob2019optimal, cornara1999satellite}.

The second echelon, the \textit{parking spares}, responds to the demands of the in-orbit spares. The parking orbits shared by all the constellations hold the inventory of parking spares. These parking orbits have a Walker Delta pattern with an inclination angle shared among the constellations. The decision maker(s) determines the number of parking orbits ($N_{\text{parking}}$) and their altitudes ($h_{\text{parking}}$), with the RAANs evenly distributed at intervals of $2\pi/N_{\text{parking}}$.  Notably, the parking orbit altitude is lower than any of the constellation altitudes ($h_{\text{parking}}<h_{\text{plane},j}\ \forall j$).

This altitude difference results in a variation in the relative RAAN drift between a parking orbit and a target orbit. The gap in RAAN shift between a parking orbit and an orbital plane of constellation $j$ is represented by the following equation, derived from Eq. \eqref{eqn: nodal precession} as
\begin{equation}\label{eqn: nodal precession gap}
    \dot{\Omega}_{\text{rel},j}=\dot{\Omega}\left(h_{\text{parking}},0,i\right)-\dot{\Omega}\left(h_{\text{plane},j},0,i\right),\quad j=1,2,\dots,m.
\end{equation}
This difference in RAAN shifts results in a periodic alignment between the parking orbit and every target orbit. A replenishment opportunity occurs whenever a parking orbit aligns with a target orbit, provided the corresponding parking orbit has available stock.

The final echelon, \textit{ground spares}, replenishes the parking spares. Ground spares are stored in warehouses or produced on demand. We assume the ground spares are always available, while the specifics of their storage and production are beyond the scope of this work. Ground spares for different satellite constellations are loaded onto a launch vehicle through ridesharing. This process is managed via a joint replenishment policy, $\left(U,\mathbf{S}\right)$, that accounts for both the total volume occupied by the satellites and the capacity of the launch vehicle.

We construct mathematical models for inventory parameters and cost elements within these multi-echelon supply-demand relationships based on the following foundational assumptions.
\begin{enumerate}
    \item In-orbit satellite failures of each satellite constellation follow a Poisson process with a constant failure rate ($\lambda_{\text{sat},j}$).
    \item Loss of spare satellites (i.e., failures of spare satellites) is not considered.
    \item In-plane spares are considered immediately available to replace a failed satellite \cite{jakob2019optimal}.
    \item A constellation's orbital plane is replenished from the spares in the closest available parking orbit (i.e., the non-empty parking spare with the earliest possible alignment). When an orbital plane's stock levels reach its ROP, a batch of spares is supplied from the closest available parking orbit \cite{jakob2019optimal}.
    \item Supply from the ground can be delivered only to a unique parking orbit \cite{jakob2019optimal}.
    \item To facilitate the tracking of the orders, an order can be processed only when no previous order is already in transit in case of in-plane spares \cite{jakob2019optimal}.
    \item Stock-outs happen rarely and thus are negligible for calculating the stock levels of in-plane spares \cite{jakob2019optimal}.
    \item Since spares have to be transferred by batches from the Earth's ground to the parking orbits and from the parking orbits to the orbital planes, the order quantity and ROP at parking orbits are assumed to be multiples of the batch size of constellation $j$ ($Q_j$) \cite{jakob2019optimal}. Thus, the values of the OUL ($S_j$) and the current inventory position ($w_j$) of a parking orbit are represented as batches of each constellation $j$.
    \item The SROP of the parking spares ($U$) is no larger than the capacity of the launch vehicle ($A$).
    \item The shipping size of a batch ($v_jQ_j$) of every constellation is less than the SROP ($U$) to ensure the feasibility of sharing launch opportunities.
    \item The inclination of the satellite constellations is identical since the change in inclination using the satellite's propulsion system requires too large $\Delta V$ for the joint replenishment strategy to be efficient.
\end{enumerate}

\subsection{Inventory Model for In-Plane Spares}
This subsection describes the details of the model for the in-plane spares, the first echelon of the proposed supply chain. The development of this model is closely tied to the concepts introduced in Section \ref{section: (s, Q) policy}.

\subsubsection{Demand Model: In-Orbit Failures}
The in-plane spares directly respond to in-orbit failures of the operational satellites. Every failure of a satellite forming the constellation is replaced by a new one waiting as an in-plane spare in the same orbital plane. The failure rate within an orbital plane of constellation $j$ ($\lambda_{\text{plane},j}$) is defined as 
\begin{equation}\label{eqn: in-planeplane failure rate}
    \lambda_{\text{plane},j} = \frac{N_{\text{sat},j} \lambda_{\text{sat},j}}{N_t}, \quad j=1,2,\dots,m,
\end{equation}
where $\lambda_{\text{sat},j}$ represents the failure rate of a satellite in constellation $j$ per year, and $N_t$ denotes the number of discretized time units per year (e.g., 365 for a daily basis, 52 for a weekly basis) \cite{jakob2019optimal}. Given a lead time $\tau$, the number of failures in an orbital plane of constellation $j$ adheres to the Poisson distribution with the following PMF:
\begin{equation}\label{eqn: in-plane demand}
    P_{\text{plane},j}(\delta;\tau)=\frac{\left(\lambda_{\text{plane},j} \tau\right)^\delta}{\delta!} e^{-\lambda_{\text{plane},j} \tau},\quad \delta\in\mathbb{Z}_+,\ j=1,2,\dots,m.
\end{equation}

\subsubsection{Replenishment Strategy: \texorpdfstring{$\left(s,Q\right)$}{} Policy with \texorpdfstring{$J_2$}{} Perturbation and Satellite's Propulsion System}

The in-orbit spares in each orbital plane of constellation $j$ are replenished using an $\left(s, Q\right)$ policy with parameters $s = s_j$ and $Q = Q_j$. Under this policy, the stock level of each orbital plane is continuously monitored, and a batch of $Q_j$ spare satellites is ordered from the nearest available parking orbit when the stock level reaches $s_j$. The time for transferring a batch of $Q_j$ satellites includes both the time required to align the orbital plane with the nearest parking orbit and the time of flight ($TOF$, in Eq. \eqref{eqn: time of flight}) for the transfer maneuver. This lead time varies depending on the relative positions of the parking orbits and the orbital plane when the order is placed.

Since the parking orbits' RAANs are equally distributed, the RAAN difference between a target orbit and its $k$th closest parking orbit {color{purple}($\Delta \Omega_k$) falls within the range
\begin{equation}\label{eqn: delta omega}
    \Delta \Omega_k \in \left[\left(k-1\right) \frac{2\pi}{N_{\text{parking}}},k \frac{2\pi}{N_{\text{parking}}}\right),\quad k=1,2,\dots,N_{\text{parking}}.
\end{equation}}
Given that the relative RAAN shift between an orbital plane of constellation $j$ and a parking orbit ($\dot{\Omega}_{\text{rel},j}$) has a relation of Eq. \eqref{eqn: nodal precession gap}, the admissible time interval for the alignment between the $k$th closest parking orbit and a target orbit of constellation $j$ ($\mathcal{T}_{j,n}$) is expressed as
\begin{equation}\label{eqn: domain of supplement time in-plane spare}
    \mathcal{T}_{j,k}=\left[\left(k-1\right) \frac{2\pi}{N_{\text{parking}}} \frac{1}{\left|\dot{\Omega}_{\text{rel},j}\right|}+t_{\text{trans},j},k \frac{2\pi}{N_{\text{parking}}} \frac{1}{\left|\dot{\Omega}_{\text{rel},j}\right|}+t_{\text{trans},j}\right),\quad \ j=1,2,\dots,m,\ k=1,2,\dots,N_{\text{parking}}
\end{equation}
where $t_{\text{trans},j}$ represents the $ToF$ for transferring from a parking orbit to an orbital plane of constellation $j$, as calculated using Eq. \eqref{eqn: time of flight}. Since the demand process and RAAN drift are independent, the alignment time is uniformly distributed within this range.

In addition, the probability that the supply comes from the $k$th closest parking orbit ($P_{\text{av},j}\left(k\right)$) is given as
\begin{equation}\label{eqn: nth parking availability}
    P_{\text{av},j}\left(k\right) = \left(1-\rho_{\text{parking}}\right)^{k-1}  \rho_{\text{parking}},\quad j=1,2,\dots,m,\ k=1,2,\dots,N_{\text{parking}},
\end{equation}
where $\rho_{\text{parking,j}}$ is the parking spares' order fill rate of constellation $j$ to be defined later.

Consequently, the probability density function (PDF) of the lead time is expressed as follows:
\begin{equation}\label{eqn: lead time distribution in-plane spare}
    g_{\text{plane},j}\left(\tau\right)=P'_{av,j}\left(k\right)  \frac{N_{\text{parking}}  \left|\dot{\Omega}_{\text{rel},j}\right|}{2\pi},\quad \tau\in \mathcal{T}_{j,n},\ j=1,2,\dots,m,\ k=1,2,\dots,N_{\text{parking}}.
\end{equation}
The right-hand side consists of two components. The first component, $P'_{\text{av},j}\left(k\right)$, is the normalized probability of supply from the $k$th closest parking orbit (i.e., $P_{\text{av},j}\left(k\right)/\sum_{k=1}^{N_{\text{parking}}}P_{\text{av},j}\left(k\right)$), and the second component is the PDF of the uniform distribution over $\mathcal{T}_{j,k}$.

\subsubsection{Inventory Parameters}
We track four key inventory parameters: order frequency, expected shortage, order fill rate, and mean stock level. The order frequency represents the number of reorders anticipated per unit time, which is derived from Eq. \eqref{eqn: order frequency} as follows:
\begin{equation}\label{eqn: order frequency plane}
    n_{\text{plane},j} = \frac{\lambda_{\text{plane},j}}{Q_j},\quad j = 1,2,\dots,m.
\end{equation}

Next, the expected shortage is obtained from Eq. \eqref{eqn: expected shortage expansion with Poisson} as
\begin{equation}\label{eqn: expected shortage in-plane}
    ES_{\text{plane},j} = \int_{\tau\ge t_{\text{trans},j}}\left(\lambda_{\text{plane},j} \tau \left( 1 - F_{\text{plane},j}(s_j - 1;\tau) \right) - s_j \left( 1 - F_{\text{plane},j}(s_j;\tau) \right)\right) g_{\text{plane},j}\left(\tau\right) \, d\tau,\quad j=1,2,\dots,m,
\end{equation}
where $F_{\text{plane},j}$ is the CDF of the Poisson distribution with a PMF of $P_{\text{plane},j}$. Following the expected shortage, the order fill rate of in-plane spares ($\rho_{\text{plane},j}$) is obtained as the fraction of the fulfilled demands by an order
\begin{equation}\label{eqn: order fill rate in-plane}
\rho_{\text{plane},j} = 1 - \frac{ES_{\text{plane},j}}{Q_j},\quad j=1,2,\dots,m.
\end{equation}

Lastly, the following represents constellation $j$'s mean stock level of in-plane spares ($\overline{SL_{\text{plane},j}}$)
\begin{equation}\label{eqn: mean stock level in-plane spare}
    \overline{SL_{\text{plane},j}} = \int_{\tau\ge t_{\text{trans},j}} \left(s_j - \lambda_{\text{plane},j} \tau + \frac{Q_j}{2} + \frac{1}{2}\right) g_{\text{plane},j}\left(\tau\right)\, d\tau,\quad j=1,2,\dots,m,
\end{equation}
derived from Eq. \eqref{eqn: mean stock level Poisson}.

\subsection{Inventory Model for Parking Spares}
This subsection focuses on the inventory model for parking spares. We build the model based on the concepts outlined in Section \ref{section: (U, S) policy}.

\subsubsection{Demand Model: Collective Demands from In-Orbit Spares}
The parking spares replenish in-orbit spares. The failure process for each orbital plane in constellation $j$ follows a Poisson process, with an order placed after every $Q_j$ failures. The time between two consecutive orders from this plane follows an Erlang-$Q_j$ distribution, with a mean of $Q_j/\lambda_{\text{plane},j}$. The parking orbits collectively replenish the demands from all orbital planes of constellation $j$. Empirically, if $N_{\text{plane},j}$ is sufficiently large (e.g., no smaller than 20), a Poisson process provides a good approximation of the stochastic demand distribution \cite{kim1991modeling,ganeshan1999managing}. Consequently, the failure rate of this Poisson process ($\lambda_{\text{parking},j}$) is given as follows:
\begin{equation}\label{eqn: parking demand rate}
    \lambda_{\text{parking},j} = \frac{N_{\text{plane},j} \lambda_{\text{plane},j}}{N_{\text{parking}} Q_j},\quad j=1,2,\dots,m.
\end{equation}
Similar to Eq. \eqref{eqn: in-plane demand}, the number of demand arrivals of constellation $j$ at a parking orbit is given as the following PMF:
\begin{equation}\label{eqn: parking demand}
    P_{\text{parking},j}(\delta;\tau)=\frac{\left(\lambda_{\text{parking},j} \tau\right)^\delta}{\delta!} e^{-\lambda_{\text{parking},j} \tau},\quad \delta\in\mathbb{Z}_+,\ j=1,2,\dots,m.
\end{equation}
Additionally, we define the total demand rate in a parking orbit ($\lambda_{\text{parking}}$) as follows:
\begin{equation}\label{eqn: total parking demand rate}
    \lambda_{\text{parking}}=\sum_{j=1}^{m}\lambda_{\text{parking},j}.
\end{equation}

\subsubsection{Replenishment Strategy: \texorpdfstring{$\left(U,\mathbf{S}\right)$}{} Policy with Launch Vehicle}
Each parking orbit is replenished by an $\left(U,\mathbf{S}\right)$ policy. In this work, the parking spares are managed in units of batches. Consequently, vector $\mathbf{a}$ represents the shipping size of each batch defined as
\begin{equation}\label{eqn: batch shipping size vector}
    \mathbf{a}=\left[v_1 Q_1,v_2 Q_2\dots,v_{m} Q_{m}\right]^\top,
\end{equation}
where $v_j$ is the shipping size of an individual satellite forming the constellation $j$, and $Q_j$ is the batch order size of constellation $j$ for in-plane spares.

The lead time for orders consists of a fixed term, $t_{\text{lau}}$, representing the time required for order processing, and a variable term accounting for the uncertainty in the timing of the next launch---an aspect not considered in the original $\left(U,\mathbf{S}\right)$ policy \cite{li2020stochastic}. This variable term follows an exponential distribution with parameter $\mu_{\text{lau}}$ \cite{jakob2019optimal}. As a result, the PDF of the lead time is expressed as follows:
\begin{equation}\label{eqn: lead time distribution parking spare}
    g_{\text{parking}}\left(\tau\right) = \frac{1}{\mu_{\text{lau}}}e^{-\left(\tau-t_{\text{lau}}\right)/\mu_{\text{lau}}},\quad \tau\in\left[t_{\text{lau}},\infty\right).
\end{equation}

\subsubsection{Inventory Parameters}
We use the MC to model and analyze the order frequency, expected order quantity, expected shortage, mean stock level, and order fill rate.

First, we define the state vector $\mathbf{w}$ based on Eq. \eqref{eqn: state definition} as follows:
\begin{equation}\label{eqn: parking state space}
    \mathcal{S}=\begin{cases}
        \mathcal{S}_0-\left\{\mathbf{0}\right\},&\text{if }\mathcal{I}_{\mathbf{w}}=\emptyset, \forall \mathbf{w}\in \mathcal{S}_0,\\
        \mathcal{S}_0,&\text{otherwise},
    \end{cases}
\end{equation}
where $\mathcal{I}_{\mathbf{w}}$ is defined in Eq. \eqref{eqn: I_w definition}. We also define the transition probability as in Eq. \eqref{eqn: state transition probability}, with $\lambda_{\text{parking},j}$ and $\lambda_{\text{parking}}$ substituting $\lambda_j$ and $\lambda$, respectively.

The MC defined by this state space and transition probability is irreducible and has a finite state space. Consequently, the vector $\mathbf{w}_\infty$, as defined by Eqs. \eqref{eqn: stationary distribution condition 1} and \eqref{eqn: stationary distribution condition 2}, represents the stationary distribution when the MC is aperiodic, or at least the limiting frequency in cases of periodicity \cite{rosenthal2006first}.

Additionally, we define the function $w_{\infty,j}\left(\cdot\right)$ that maps the inventory state of constellation $j$'s parking spares ($w_j$) to its probability in the stationary distribution as
\begin{equation}\label{eqn: marginalized stationary distribution}
    w_{\infty,j}\left(\omega\right)=\sum_{\mathbf{w}\in\mathcal{S},w_j=\omega}w_{\infty,\mathbf{w}},\quad \omega\in\mathbb{Z}_+,\ j=1,2,\dots,m.
\end{equation}

In the infinite time horizon, the outflow and inflow are balanced for every state. Therefore, the outflows and inflows from states $\mathbf{e}_j$ for $j=1,2,\dots,m$ are equal, and these represent the order frequency. Consequently, the order frequency ($n_{\text{parking}}$) is given as
\begin{equation}\label{eqn: order frequency parking}
    n_{\text{parking}}=\lambda_{\text{parking}}\sum_{j=1}^mw_{\infty, \mathbf{e}_j}.
\end{equation}

Unlike in-plane spares, the order quantity of parking spares is not predetermined but is determined when an order is triggered. To calculate its expected value, we define $\mathcal{I}_{j,1}$ and $\mathcal{I}_{j,2}$ as
\begin{equation}\label{eqn: I_j1 definition}
    \mathcal{I}_{j,1}=\left\{\mathbf{w}\in \mathcal{S}:\mathbf{w}+\mathbf{e}_j\notin \mathcal{S},\ \mathbf{a}^\top \left(\mathbf{w}+\mathbf{e}_j\right)\le A\right\},\quad j=1,2,\dots,m,
\end{equation}
\begin{equation}\label{eqn: I_j2 definition}
    \mathcal{I}_{j,2}=\left\{\mathbf{w}\in \mathcal{S}:\mathbf{w}+\mathbf{e}_j\notin \mathcal{S},\ \mathbf{a}^\top \left(\mathbf{w}+\mathbf{e}_j\right)> A\right\},\quad j=1,2,\dots,m.
\end{equation}
These sets represent the states where a demand from the in-plane spares of constellation $j$ triggers an order. $\mathcal{I}_{j,1}$ corresponds to states where all cumulative demands can be fulfilled by a single launch. In contrast, $\mathcal{I}_{j,2}$ represents states where the cumulative demands cannot be fully fulfilled by a single launch, but the last demand must be included in the next launch. Therefore, if the current state $\mathbf{w}$ is in $\mathcal{I}_{j,1}$, the order quantity corresponds to the vector $\mathbf{w} + \mathbf{e}_j$; if in $\mathcal{I}_{j,2}$, the order quantity follows $\mathbf{w}$. As a result, we formulate the expected order quantity ($EQ_{\text{parking},j}$) as
\begin{equation}\label{eqn: expected order quantity}
    EQ_{\text{parking},j} = \sum_{j'=1}^m\left(\lambda_{\text{parking},j'}\mathbf{e}_j^\top \left(\sum_{\mathbf{w}\in \mathcal{I}_{j',1}}\left(\mathbf{w}+\mathbf{e}_{j'}\right) + \sum_{\mathbf{w}\in \mathcal{I}_{j',2}}\mathbf{w}\right)\right),\quad j=1,2,\dots,m.
\end{equation}

Next, we introduce additional relations to express the expected shortage and mean stock level. Given a lead time, $\tau$, the demands during this lead time, $\delta$, satisfy
\begin{equation}\label{eqn: net stock additional relation 1}
    NS_j\left(t+\tau\right)=NS_j\left(t\right)+O_j\left(t\right)-\delta,
\end{equation}
where $t$ represents the current time, $NS_j(t)$ and $NS_j(t+\tau)$ are the net stocks at time $t$ and $t+\tau$, respectively, and $O_j(t)$ is the on-order stock at time $t$ \cite{zipkin1986stochastic}. Since Eqs. \eqref{eqn: net stock definition} and \eqref{eqn: state definition} hold, the following relation also holds:
\begin{equation}\label{eqn: net stock additional relation 2}
    NS_j\left(t+\tau\right) = S_j-w_j\left(t\right)-\delta,
\end{equation}
where $w_j(t)$ is the $j^{\text{th}}$ component of the state vector at time $t$. Since the stationary distribution exists, as $t \rightarrow \infty$, the following relation holds \cite{zipkin1986stochastic}:
\begin{equation}\label{eqn: net stock additional relation 3}
    NS_{\infty,j}\left(\omega\right)=S_j-w_{\infty,j}\left(\omega\right)-\delta,
\end{equation}
where $w_{\infty,j}\left(\cdot\right)$ is the function defined in Eq. \eqref{eqn: marginalized stationary distribution}, and $NS_{\infty,j}\left(\cdot\right)$ is the function that maps the net stock level to its probability at the stationary distribution.

Based on Eqs. \eqref{eqn: net stock additional relation 3} and \eqref{eqn: net stock additional relations}, we formulate the expected shortage ($ES_{\text{parking},j}$) and the mean stock level ($\overline{SL_{\text{parking},j}}$) of constellation $j$'s parking spares as follows \cite{li2020stochastic}:
\begin{equation}\label{eqn: expected shortage parking}
    ES_{\text{parking},j} = \int_{\tau\ge t_{\text{lau}}}\sum_{\omega\in\mathbb{Z}_+}\sum_{\delta\in\mathbb{Z}_+}\left(S_j-\omega-\delta\right)^- w_{\infty,j}\left(\omega\right)P_{\text{parking},j}(\delta;\tau) g_{\text{parking}}(\tau)\, d\tau,
\end{equation}
\begin{equation}\label{eqn: mean stock level parking}
    \overline{SL_{\text{parking},j}} = S_j-\sum_{\omega\in\mathbb{Z}_+} \omega w_{\infty,j}\left(\omega\right) - \int_{\tau\ge t_{\text{lau}}}\sum_{\delta\in\mathbb{Z}_+}\delta P_{\text{parking},j}(\delta;\tau) g_{\text{parking}}(\tau)\, d\tau + ES_{\text{parking},j},
\end{equation}
where $\left(S_j-\omega-\delta\right)^-=\max\left\{-S_j+\omega+\delta,0\right\}$. Following these relations, we introduce the order fill rate for constellation $j$ ($\rho_{\text{parking},j}$) as follows:
\begin{equation}\label{eqn: order fill rate parking}
\rho_{\text{parking},j}=1-\frac{ES_{\text{parking},j}}{EQ_{\text{parking},j}}.
\end{equation}

\subsection{Total Maintenance Cost}
The elements of maintenance cost are modeled using parameters introduced in the previous section. The cost components include launch, holding, maneuvering, and manufacturing costs \cite{jakob2019optimal}. Since $m$ constellations share launch opportunities, the launch costs must also be distributed among the constellations.

First, we define the annual maintenance cost of constellation $j$ except the launch cost ($C_j$) as
\begin{equation}\label{eqn: individual cost}
    C_j = C_{\text{hold},j} + C_{\text{maneuv},j} + C_{\text{manufac},j},\quad j=1,2,\dots,m,
\end{equation}
where the annual holding ($C_{\text{hold},j}$), maneuvering ($C_{\text{maneuv}}$), and manufacturing ($C_{\text{manufac}}$) costs are calculated as follows:
\begin{equation}\label{eqn: holding cost}
    C_{\text{hold},j} = c_{\text{hold},j} \left(\overline{SL_{\text{plane},j}}{N_{\text{plane},j}} + \overline{SL_{\text{parking},j}} Q_j{N_{\text{parking}}}\right),\quad j=1,2,\dots,m,
\end{equation}
\begin{equation}\label{eqn: maneuvering cost}
    C_{\text{maneuv},j} = M_{\text{fuel},j} \epsilon_{\text{maneuv},j} n_{\text{plane},j} N_{\text{plane},j} Q_j N_t,\quad j=1,2,\dots,m,
\end{equation}
\begin{equation}\label{eqn: manufacturing cost}
    C_{\text{manufac},j} = c_{\text{manufac},j}\lambda_{\text{sat},j} N_{\text{plane},j} N_{\text{sat},j} ,\quad j=1,2,\dots,m.
\end{equation}
Here, the holding cost is calculated as the product of the annual inventory holding cost per unit ($c_{\text{hold},j}$) and the sum of the expected inventory levels of in-plane spares and parking spares ($\overline{SL_{\text{plane},j}} + \overline{SL_{\text{parking},j}} Q_j$). The maneuvering cost is determined by multiplying the cost of maneuvering per event ($M_{\text{fuel},j} \epsilon_{\text{maneuv},j}$) by the total number of annual maneuvering events ($n_{\text{plane},j} N_{\text{plane},j} Q_j N_t$). Finally, the manufacturing cost is the product of the manufacturing cost per satellite ($c_{\text{manufac},j}$) and the total number of annual satellite failures ($\lambda_{\text{sat},j} N_{\text{plane},j} N_{\text{sat},j}$).

In addition to these cost components, we introduce $y_j$ ($\ge0$), $j = 1, 2, \dots, m$, representing the fraction of the launch costs allocated to constellation $j$, where the sum of all $y_j$ equals 1:
\begin{equation}\label{eqn: launch cost fraction}
    \sum_{j=1}^m y_j = 1.
\end{equation}
Based on this, we define the total expected spare strategy annual cost (TESSAC \cite{jakob2019optimal}) for constellation $j$ as $TESSAC_j$, and the total TESSAC for all constellations as $TESSAC$, as follows:
\begin{equation}\label{eqn: TESSAC j}
    TESSAC_j = y_j C_{\text{lau}} + C_j, \quad j = 1,2,\dots,m,
\end{equation}
\begin{equation}\label{eqn: TESSAC}
    TESSAC = \sum_{j=1}^m TESSAC_j = C_{\text{lau}} + \sum_{j=1}^m C_j,
\end{equation}
where we define the annual launch cost ($C_{\text{lau}}$) as the product of the launch cost ($c_{\text{lau}}$) and the annual expected number of launches ($n_{\text{parking}} N_{\text{parking}} N_t$):
\begin{equation}\label{eqn: launch cost}
    C_{\text{lau}} = c_{\text{lau}} n_{\text{parking}} {N_{\text{parking}}} N_t.
\end{equation}

\subsection{Model Validation}
We validate the mathematical model of inventory parameters and relevant costs by comparing the model's results to those of the simulation. For this purpose, we generate 25 random instances for each value of $m$, varying $m$ from 2 to 5. Tables \ref{tab: common parameters} and \ref{tab: sampled trade space} summarize the common and varying parameters, respectively.

\begin{table}[hbt!]
    \centering
    \caption{Common parameters of instances generated for validation study, $j=1,2,\dots,m$}
    \begin{tabular}{lccc}
        \hline\hline
        Parameter & Notation & Value & Unit\\\hline
        Fuel mass conversion coefficient & $\epsilon_{\text{fuel},j}$ & 0.01 & M\$/kg\\
        Annual satellite holding cost & $c_{\text{hold},j}$ & 0.5 & M\$/satellite/year\\
        Exhaust velocity of satellite's propulsion system & $V_{\text{ex},j}$ & 11.77 & km/s\\
        Specific impulse of satellite's propulsion system & $I_{\text{sp},j}$ & 1200 & sec\\
        Launch cost & $c_{\text{lau}}$ & 200 & M\$\\
        Launch capacity & $A$ & 250 & \# slots\\
        Dry mass coefficient & $\alpha_1$ & 150 & kg/slot\\
        Propellant mass rate coefficient& $\alpha_2$ & 1.3E-3 & kg/slot/s\\
        Production cost coefficient& $\alpha_3$ & 0.5 & \$M/slot\\
        \hline\hline
    \end{tabular}
    \label{tab: common parameters}
\end{table}

\begin{table}[hbt!]
    \centering
    \caption{Sampled trade space of varying parameters for validation scenarios, $j=1,2,\dots,m$}
    \begin{tabular}{lccc}
        \hline\hline
        Parameter & Notation & Range & Unit\\\hline
        Launch order processing time & $t_{\text{lau}}$ & $\left[20, 80\right]$ & weeks\\
        Mean waiting time to launch & $\mu_{\text{lau}}$ & $\left[20, 80\right]$ & weeks\\
        Inclination & $i$ & $\left[40, 80\right]$ & deg\\
        Satellite failure rate & $\lambda_{\text{sat},j}$ & $\left[0.05, 0.2\right]$ & failures/year\\
        Number of orbital planes & $N_{\text{plane},j}$ & $\left[20, 40\right]$ & planes\\
        Number of parking orbits & $N_{\text{parking}}$ & $\left[1, 20\right]$ & planes\\
        Number of satellites per orbital plane & $N_{\text{sat},j}$ & $\left[20, 60\right]$ & satellites/plane\\
        Orbital plane altitude & $h_{\text{plane},j}$ & $\left[500, 2000\right]$ & km\\
        Parking orbit altitude & $h_{\text{parking}}$ & $\left[400, 1000\right]$ & km\\
        ROP of in-plane spares & $s_j$ & $\left[1, 20\right]$ & satellites\\
        Order quantity of in-plane spares & $Q_j$ & $\left[1, 40\right]$ & satellites\\
        OUL of parking spares & $S_j$ & $\left[1, 40\right]$ & - \\
        SROP of parking spares & $U$ & $\left[10, 250\right]$ & \# slots\\
        Shipping size of individual satellite & $v_j$ & $\left[1, 4\right]$ & \# slots\\
        \hline\hline
    \end{tabular}
    \label{tab: sampled trade space}
\end{table}

For each value of $m$, we sample parameters uniformly at random until 25 instances satisfying Eqs. \eqref{eqn: condition 1}-\eqref{eqn: condition 5} are generated:
\begin{equation}\label{eqn: condition 1}
    s_j - Q_j \le 0,\quad j=1,2,\dots,m,
\end{equation}
\begin{equation}\label{eqn: condition 2}
    \sum_{j=1}^mS_jv_jQ_j\ge U,
\end{equation}
\begin{equation}\label{eqn: condition 3}
    v_j Q_j+1\le U,\quad j=1,2,\dots,m,
\end{equation}
\begin{equation}\label{eqn: condition 4}
    \rho_{\text{parking},j}\ge \rho_{\text{parking},j}^{\text{req}},\quad j=1,2,\dots,m,
\end{equation}
\begin{equation}\label{eqn: condition 5}
    \rho_{\text{plane},j}\ge \rho_{\text{plane},j}^{\text{req}},\quad j=1,2,\dots,m,
\end{equation}
These conditions are rooted in the model assumptions. Specifically, Eqs. \eqref{eqn: condition 1} and \eqref{eqn: condition 2} are for ensuring the cyclic behavior of the inventory management system. Eq. \eqref{eqn: condition 3} allows for the joint replenishment strategy to be feasible, while Eqs. \eqref{eqn: condition 4} and \eqref{eqn: condition 5} enforce rare shortages. The parameters $\rho_{\text{parking},j}^{\text{req}}$ and $\rho_{\text{plane},j}^{\text{req}}$ denote the thresholds for the order fill rate, which are set at 0.98 for all $j$ in the validation study.

In addition, we assume that the dry mass ($M_{\text{dry},j}$), propellant mass rate ($\dot{M}_{\text{prop,j}}$), and production cost ($c_{\text{sat},j}$) of a satellite are proportional to its shipping size ($v_j$) expressed as follows:
\begin{equation}\label{eqn: dry mass relation}
    M_{\text{dry},j} = \alpha_1 v_j,\quad j=1,2,\dots,m,
\end{equation}
\begin{equation}\label{eqn: mdot relation}
    \dot{M}_{\text{prop},j} = \alpha_2 v_j,\quad j=1,2,\dots,m,
\end{equation}
\begin{equation}\label{eqn: production cost relation}
    c_{\text{sat},j} = \alpha_3 v_j,\quad j=1,2,\dots,m,
\end{equation}
where the values of $\alpha_1$, $\alpha_2$, and $\alpha_3$ are presented in Table \ref{tab: common parameters}.

We compare the model's results with those obtained from 100 Monte Carlo simulations for each scenario and calculate the corresponding errors. In the inventory simulation, stochastic elements---failures and lead time---are incorporated. The relative error, which evaluates the order frequency and TESSAC, is determined using the formula in Eq. \eqref{eqn: model value error} \cite{ganeshan1999managing}:
\begin{equation}\label{eqn: model value error}
    \frac{\left| value_{\text{sim}} - value_{\text{model}} \right|}{value_{\text{sim}}} \times 100.
\end{equation}

In contrast, the demand rate of parking spares and the mean stock level of each spare are calculated using
\begin{equation}\label{eqn: model value error max}
    \underset{j=1,2,\dots,m}{\text{max}} \frac{\left| value_{\text{sim},j} - value_{\text{model},j} \right|}{value_{\text{sim},j}} \times 100.
\end{equation}
With this relation, the maximum errors among the constellations are selected for the average error calculation. That is, in an instance of $m=3$, if the errors of $\overline{SL_{\text{plane},j}}$ are 0.3\%, 0.6\%, and 0.4\% for $j=1,2,3$, respectively, then only 0.6\% is accounted for when calculating the mean value. Consequently, the validation result reflects the worst errors among different constellations in each instance, setting a more challenging standard for the accuracy test.

Similarly, rational parameters, such as the order fill rate, are evaluated using the absolute error, which is calculated as
\begin{equation}\label{eqn: model rate error max}
\underset{j=1,2,\dots,m}{\text{max}}\left|value_{\text{sim},j} - value_{\text{model},j}\right|.
\end{equation}

The validation results presented in Table \ref{tab: validation results} indicate that, although the errors tend to increase with higher values of $m$, the model maintains sufficient accuracy, with relative errors below 3\% and absolute errors below 0.5\%.
\begin{table}[hbt!]
    \centering
    \caption{Averaged errors of result parameters, $j=1,2,\dots,m$}
    \begin{tabular}{lccccc}
        \hline\hline
        Result Parameter & Notation & \multicolumn{4}{c}{Mean Maximum Error} \\
         & & $m=2$ & $m=3$ & $m=4$ & $m=5$ \\\hline
        \textbf{Relative Errors (\%)} & & & & & \\
        Demand rate of parking spares & $\lambda_{\text{parking},j}$ & 1.19 & 1.44 & 1.74 & 1.56 \\
        Mean stock level of in-plane spares & $\overline{SL_{\text{plane},j}}$ & 0.43 & 0.51 & 0.48 & 0.50 \\
        Mean stock level of parking spares & $\overline{SL_{\text{parking},j}}$ & 0.60 & 0.63 & 1.29 & 0.93 \\
        Order frequency of parking spares & $n_\text{parking}$ & 1.10 & 1.63 & 1.92 & 2.26 \\
        Total relevant cost & $TESSAC$ & 0.17 & 0.21 & 0.25 & 0.20 \\
        \textbf{Absolute Errors (\%p)} & & & & & \\
        Order fill rate of in-plane spares & $\rho_{\text{plane},j}$ & 0.01 & 0.02 & 0.02 & 0.03 \\
        Order fill rate of parking spares & $\rho_{\text{parking},j}$ & 0.05 & 0.08 & 0.36 & 0.12 \\
        \hline\hline
    \end{tabular}
    \label{tab: validation results}
\end{table}

\section{Decision-Making Scenarios and Solution Procedures}\label{sec: section 4}
We introduce two decision-making scenarios: centralized and decentralized. Each scenario represents a distinct decision-making context and, accordingly, requires a unique solution procedure.

\subsection{Centralized Scenario}
Imagine a situation where a single authority maintains $m$ different constellations. For instance, a single company might operate multiple constellations with varying characteristics, or multiple companies might own individual constellations while delegating system maintenance to an external authority. This authority would have full access to information and control over all maintenance-related activities. The latter is common in terrestrial logistics, where third-party logistics providers (3PLs) play a significant role \cite{chopra2019supply}, and it is becoming increasingly relevant. A few studies have also begun to address the role of such players in the market \cite{ho2024space, hwang2024optimal}. In this centralized context, maintenance decisions can be made collectively to maximize the efficiency of the whole supply chain.

We formulate the decision-making problem in this centralized scenario as the optimization problem $\mathbf{P}_{\text{JC}}$, where the objective is to minimize $TESSAC$ while maintaining the required service level for spares in each constellation. We use the term ``strategy'' to refer to decisions regarding parking orbit configurations and replenishment policies:
\newline\newline\noindent
($\mathbf{P}_{\text{JC}}$) Joint Replenishment Strategy Optimization in a Centralized Scenario
\begin{equation}\label{eqn: JC obj}
    TESSAC\left(\mathbf{x}_{\text{JC}}^*;\mathbf{p}\right)=\underset{\mathbf{x}}{\text{min}}\ TESSAC\left(\mathbf{x};\mathbf{p}\right)
\end{equation}
\noindent
subject to
\begin{center}
    Eqs. \eqref{eqn: condition 1}--\eqref{eqn: condition 5},
\end{center}
\begin{equation}\label{eqn: JC const 1}
    \mathbf{x}=[\mathbf{x}_{\text{sh}},\mathbf{x}_1,\dots,\mathbf{x}_{m}],
\end{equation}
\begin{equation}\label{eqn: JC const 2}
    \mathbf{x}_{\text{sh}}=\left[h_{\text{parking}},N_{\text{parking}},U\right]\in \mathcal{H}_{\text{parking}}\times\mathcal{N}_{\text{parking}}\times\mathbb{Z}_{++},
\end{equation}
\begin{equation}\label{eqn: JC const 3}
    \mathbf{x}_{j}=\left[s_j,Q_j,S_j\right]\in\mathbb{Z}_{++}^{3},\quad j=1,2,\dots,m,
\end{equation}
where $\mathbf{p}$ denotes the parameters external to the replenishment strategy design but essential for modeling the system, such as the specifications of the satellite propulsion system and the failure rate of the satellites. $\mathcal{H}_{\text{parking}}$ and $\mathcal{N}_{\text{parking}}$ represent the admissible value sets for the altitude and the number of parking orbits, respectively. $\mathbf{x}_{\text{sh}}$ is the vector containing the shared decision variables, and $\mathbf{x}_j$ represents the replenishment policy parameters for constellation $j$.

\subsection{Decentralized Scenario}
The decentralized scenario assumes multiple authorities, where each authority manages its associated constellation. It seeks to minimize its own maintenance costs, making it more challenging to reconcile multiple, potentially conflicting objectives than in a centralized setting. To tackle this complexity, we propose a two-step procedure. First, efficient (Pareto-Optimal) solutions are identified as potential alternatives (Step 1). Then, the operators negotiate to select a single agreed solution from the alternatives (Step 2). During this negotiation, the bargaining power of each participant significantly influences the outcome. The remainder of this subsection discusses each step in detail.
\subsubsection{Step 1: Generating the Group of Efficient Solutions}
In the decentralized scenario, each operator aims to minimize its own maintenance cost, $TESSAC_j$ ($j=1,2,\dots,m$). These distinct objectives form the basis of the multi-objective optimization problem $\mathbf{P}_{\text{JD}}^1$, minimizing each operator's maintenance cost is treated as a separate objective function.
\newline\newline\noindent
($\mathbf{P}_{\text{JD}}^1$) Joint Replenishment Strategy Optimization in a Decentralized Scenario
\begin{equation}\label{eqn: JD obj}
    \mathcal{X}^{\text{efficient}}_{\text{JD}}=\min_{\mathbf{x},\mathbf{y}}\ \left[TESSAC_1,TESSAC_2,\dots,TESSAC_m\right]\left(\mathbf{x},\mathbf{y};\mathbf{p}\right)
\end{equation}
\noindent
subject to
\begin{center}
    Eqs. {\eqref{eqn: launch cost fraction}, \eqref{eqn: condition 1}--\eqref{eqn: condition 5}, \eqref{eqn: JC const 1}--\eqref{eqn: JC const 3},}
\end{center}
\begin{equation}\label{eqn: JD const 1}
    TESSAC_j\left(\mathbf{x},\mathbf{y};\mathbf{p}\right)\le TESSAC_{\text{ref},j},\quad j=1,2,\dots,m,
\end{equation}
\begin{equation}\label{eqn: JD const 2}
    \mathbf{y}\in \mathbb{R}_+^{m},
\end{equation}
where we introduce the launch cost fractions ($\mathbf{y}$) as decision variables, and the set of efficient solutions (i.e., Pareto optimal solutions), $\mathcal{X}^{\text{eff}}_{\text{JD}}=\left\{\left(\mathbf{x}_{\text{JD}}^{\text{eff}},\mathbf{y}_{\text{JD}}^{\text{eff}}\right)\right\}$, are obtained as the results. In addition, we introduce Eq. \eqref{eqn: JD const 1} to ensure that $TESSAC_j$ remains less than the reference value $TESSAC_{\text{ref},j}$, which is the maximum expenditure each operator is willing to incur to participate in the joint replenishment strategy. The reference cost can be set in various ways; for example, it could represent the total maintenance cost if each operator maintained their constellation individually.

A solution is considered efficient if not dominated by any other solution, meaning no feasible solution is superior to it across all operators' perspectives. That is, for two feasible solutions, $\left(\mathbf{x},\mathbf{y}\right)$ and $\left(\mathbf{x}',\mathbf{y}'\right)$, the solution $\left(\mathbf{x},\mathbf{y}\right)$ is said to be dominated by $\left(\mathbf{x}',\mathbf{y}'\right)$ if the following conditions are met:
 \begin{equation}
     TESSAC_j\left(\mathbf{x},\mathbf{y};\mathbf{p}\right) \ge TESSAC_j\left(\mathbf{x}',\mathbf{y}';\mathbf{p}\right), \quad \forall j\in\left\{1,2,\dots,m\right\},
 \end{equation}
 \begin{equation}
     TESSAC_j\left(\mathbf{x},\mathbf{y};\mathbf{p}\right) > TESSAC_j\left(\mathbf{x}',\mathbf{y}';\mathbf{p}\right), \quad \exists j\in\left\{1,2,\dots,m\right\}.
 \end{equation}
In this case, there is no reason for the operators to select $\left(\mathbf{x},\mathbf{y}\right)$ over $\left(\mathbf{x}',\mathbf{y}'\right)$ as for agreement. Consequently, it can be assumed that the solution for agreement is chosen exclusively from $\mathcal{X}_{\text{JD}}^{\text{eff}}$.

\subsubsection{Step 2: Choosing a Solution for Agreement}
In the second step, a solution is chosen from the set of efficient solutions to form the final agreement. During this phase, all operators participate in the negotiation to ensure their interests are accurately represented in the outcome. We denote the proxy for bargaining power of each operator by $\beta_j$, and propose a straightforward selection mechanism that minimizes the weighted sum of $TESSAC_j$ with $\beta_j$ ($\ge0$) formulated as the following optimization problem $\mathbf{P}_{\text{JD}}^2$. The weights satisfy $\sum_{j=1}^m\beta_j=1$.
\newline\newline\noindent
($\mathbf{P}_{\text{JD}}^2$) Solution Selection for Agreement
\begin{equation}\label{eqn: step 2}
    \left(\mathbf{x}_{\text{JD}}^*,\mathbf{y}_{\text{JD}}^*\right)=\underset{\left(\mathbf{x},\mathbf{y}\right)\in\mathcal{X}_{\text{JD}}^{\text{eff}}}{\text{argmin}}\sum_{j=1}^m \beta_j TESSAC_j\left(\mathbf{x},\mathbf{y};\mathbf{p}\right).
\end{equation}

\begin{figure}[hbt!]
\centering
\includegraphics[page=4, width=1.0\textwidth]{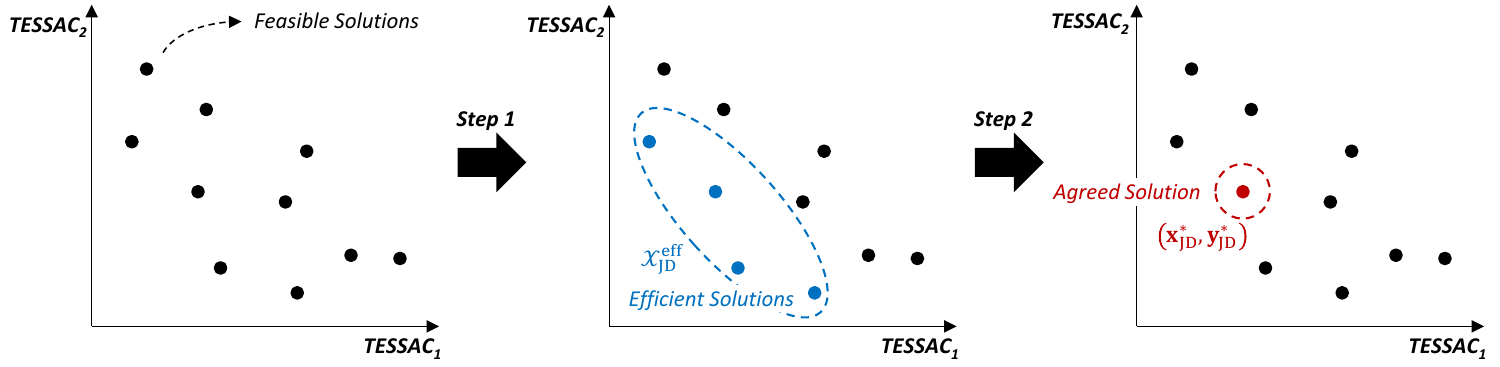}
\caption{Solution procedure for decentralized scenario, $m=2$}
\label{fig: solution procedure for decentralized context}
\end{figure}

Figure \ref{fig: solution procedure for decentralized context} summarizes the process with an example of $m=2$.

\section{Case Study}\label{sec: section 5}
We conduct two case studies to provide valuable insights and demonstrate the applicability of the proposed model and solution framework. First, we perform a parametric study to identify conditions under which the joint replenishment strategy is effective and could be considered. Second, we analyze a scenario involving the joint replenishment of three constellations to illustrate the solution procedure.

We consider two launch options for this case study: a normal launcher and a mega-launcher. The normal launcher is an analogy for existing launch vehicles such as Falcon 9, Soyuz, and Ariane 5. In contrast, the mega-launcher represents a significantly larger vehicle anticipated in the near future, such as Starship, capable of carrying greater payloads at a lower specific cost, albeit with a longer lead time. By utilizing a mega-launcher, a satellite constellation operator can benefit from reduced satellite launch costs---\textit{economies of scale}. However, this approach may require the operator to maintain a larger inventory of spare satellites due to the extended lead time and increased capacity per launch. Consequently, sharing launch opportunities might offer additional advantages in certain scenarios. The parameters for the launch options are summarized in Table \ref{tab: launch service parameters}.

\begin{table}[hbt!]
    \centering
    \caption{Parameters of the launch services for the case study}
    \begin{tabular}{lcccc}
        \hline\hline
        Parameter & Notation & \multicolumn{2}{c}{Value} & Unit\\
        & & Normal Launcher & Mega Launcher & \\\hline
        Launch cost & $c_{\text{lau}}$ & 67 & 200 & \$M\\
        Launch capacity & $A$ & 40 & 250 & \# slots\\
        Launch order processing time & $t_{\text{lau}}$ & 12 & 32 & weeks\\
        Mean waiting time to launch & $\mu_{\text{lau}}$ & 8 & 20 & weeks\\
        \hline\hline
    \end{tabular}
    \label{tab: launch service parameters}
\end{table}

We adopt the Genetic Algorithm (GA) \cite{holland1992genetic}, a metaheuristic, and NSGA-II \cite{deb2002fast}, a GA-based multi-objective metaheuristic, for solving all optimization problems in this case study.

\subsection{Case Study 1: Parametric Study}
In the parametric study, we assess how variations among different constellations influence the joint replenishment strategy. We begin by establishing a set of reference parameters and then vary one parameter in each scenario, labeled A through D. Next, we compare the resulting annual maintenance costs under joint replenishment with those under the independent replenishment strategy. The reference parameters are provided in Table \ref{tab: reference parameter 1}, and Equations \eqref{eqn: dry mass relation}--\eqref{eqn: mdot relation} apply throughout.
\begin{table}[hbt!]
    \centering
    \caption{Reference parameters for parametric study}
    \begin{tabular}{lccc}
        \hline\hline
        Parameter & Notation & Value & Unit\\\hline
        Inclination & $i$ & 60 & deg\\
        Satellite failure rate & $\lambda_{\text{sat}}$ & 0.2 & failures/year\\
        Number of orbital planes in the constellation & $N_{\text{plane}}$ & 30 & planes\\
        Number of satellites per orbital plane & $N_{\text{sat}}$ & 30 & satellites/plane\\
        Altitude of orbital planes in the constellation & $h_{\text{plane}}$ & 1200 & km\\
        Shipping size of individual satellite & $v$ & 1 & \# slots\\
        Fuel mass conversion coefficient & $\epsilon_{\text{fuel}}$ & 0.01 & \$M/kg\\
        Annual satellite holding cost & $c_{\text{hold}}$ & 0.5 & \$M/satellite/year \\
        Exhaust velocity of satellite's propulsion system & $V_{\text{ex}}$ & 11.77 & km/s\\
        Specific impulse of satellite's propulsion system & $I_{\text{sp}}$ & 1200 & sec\\
        Dry mass coefficient & $\alpha_1$ & 150 & kg/slot\\
        Propellant mass rate coefficient & $\alpha_2$ & 0.0013 & kg/slot/s\\
        Production cost coefficient & $\alpha_3$ & 0.5 & \$M/slot\\
        \hline\hline
    \end{tabular}
    \label{tab: reference parameter 1}
\end{table}

Four parameters and their corresponding cases are analyzed. Case 1-A varies the shipping size of an individual satellite ($v_2$) across values of 1, 2, 3, and 4, indexed from 0 to 3. Case 1-B varies the number of satellites per plane ($N_{\text{sat},2}$) across values of 30, 40, 50, and 60, indexed from 0 to 3. This variation leads to differing demand rates and has a similar effect to varying the failure rate of the satellites ($\lambda_{\text{sat},2}$). Case 1-C varies the altitude of the orbital planes of the constellation ($h_{\text{plane},2}$) 1200, 1300, 1400, and 1500, indexed from 0 to 3. This variation affects the lead time from the parking spares to the in-plane spares. Case 1-D varies the number of planes in the satellite constellation ($N_{\text{plane},2}$) across values of 30, 40, 50, and 60, indexed from 0 to 3. While similar to Case 1-B in influencing demand rates, it also changes the number of locations where demand occurs---an effect that cannot be fully explained by demand rates alone. Table \ref{tab: parametric study case summary} summarizes the cases along with their corresponding ground analogies.
\begin{table}[hbt!]
    \centering
    \caption{Summary of the parameters for Case Study 1}
    \begin{tabular}{cclcc}
        \hline\hline
        Case ID & Parameter & Ground Analogy & Values (index 0 to 3) & Units\\\hline
        Case 1-A & $v_2$ & shipping size of SKU & 1, 2, 3, 4 & \# slots\\
        Case 1-B & $N_{\text{sat},2}$ & demand rate at each retailer & 30, 40, 50, 60 & satellites/plane\\
        Case 1-C & $h_{\text{plane},2}$ & locations of retailers & 1200, 1300, 1400, 1500 & km\\
        Case 1-D & $N_{\text{plane},2}$ & the number of retailers & 30, 40, 50, 60 & planes\\
        \hline\hline
    \end{tabular}
    \label{tab: parametric study case summary}
\end{table}

In each instance of every case, we consider replenishing two constellations: one follows the reference parameters, while the other adopts the varying parameters. The optimal solution from the independent strategy, where in-plane spares and parking spares are replenished using $\left(s, Q\right)$ policies with parameters $\left(s_1,Q_1\right)$ and $\left(k_{R_1,3}, k_{Q_1,3}\right)$, respectively, as described by Jakob et al. \cite{jakob2019optimal}, is compared with the optimal solution from the joint replenishment strategy, obtained by solving $\mathbf{P}_{\text{JC}}$. All scenarios assume the use of a mega-launcher whose parameters are listed in Table \ref{tab: launch service parameters}.

The decision variables for the independent and joint replenishment strategies are presented in Tables \ref{tab: parametric study independent decision variables} and \ref{tab: parametric study joint decision variables}, respectively. Here, $Q_{\text{max}}=\lfloor A/v\rfloor$ and $Q_{\text{max},j}=\lfloor A/v_j\rfloor$ for $j=1,2$. The admissible values for the altitude of parking orbits are discretized in 50 km intervals within the specified range. Thus, $\mathcal{H}_{\text{parking}}$ is defined as
\begin{equation}\label{eqn: parking orbits values}
    \mathcal{H}_{\text{parking}}=\left\{500,550,\dots,1000\right\}.
\end{equation}

\begin{table}[hbt!]
    \centering
    \caption{Decision variables of independent replenishment strategy for parametric study}
    \begin{tabular}{lccc}
        \hline\hline
        Parameter & Notation & Range & Unit\\\hline
        ROP of in-plane spares & $s_1$ & $\left[1, 10\right]$ & satellites\\
        Order quantity of in-plane spares & $Q_1$ & $\left[1, Q_{\text{max}}\right]$ & satellites\\
        ROP of parking spares & $k_{R_1,3}$ & $\left[1, 40\right]$ & - \\
        Order quantity of parking spares & $k_{Q_1,3}$ & $\left[1, Q_{\text{max}}\right]$ & - \\
        Number of parking orbits & $N_{\text{parking}}$ & $\left[1, 20\right]$ & planes\\
        Altitude of parking orbits & $h_{\text{parking}}$ & $\left[500, 1000\right]$ & km\\
        \hline\hline
    \end{tabular}
    \label{tab: parametric study independent decision variables}
\end{table}

\begin{table}[hbt!]
    \centering
    \caption{Decision variables of joint replenishment strategy for parametric study, $j=1,2$}
    \begin{tabular}{lccc}
        \hline\hline
        Parameter & Notation & Range & Unit\\\hline
        ROP of in-plane spares & $s_j$ & $\left[1, 10\right]$ & satellites\\
        Order quantity of in-plane spares & $Q_j$ & $\left[1, Q_{\text{max},j}\right]$ & satellites\\
        OUL of parking spares & $S_j$ & $\left[1, 40\right]$ & - \\
        SROP of parking spares & $U$ & $\left[200, 250\right]$ & \# slots\\
        Number of parking orbits & $N_{\text{parking}}$ & $\left[1, 20\right]$ & planes\\
        Altitude of parking orbits & $h_{\text{parking}}$ & $\left[500, 1000\right]$ & km\\
        \hline\hline
    \end{tabular}
    \label{tab: parametric study joint decision variables}
\end{table}

For the independent strategy, the two constellations' TESSAC values are summed to obtain a total measure, while for the joint replenishment strategy, the TESSAC is computed according to Eq. \eqref{eqn: TESSAC}. The relative difference between these two totals is then calculated to evaluate the cost-effectiveness of joint replenishment.

\begin{figure}[hbt!]
\centering
\includegraphics[page=5, width=0.8\textwidth]{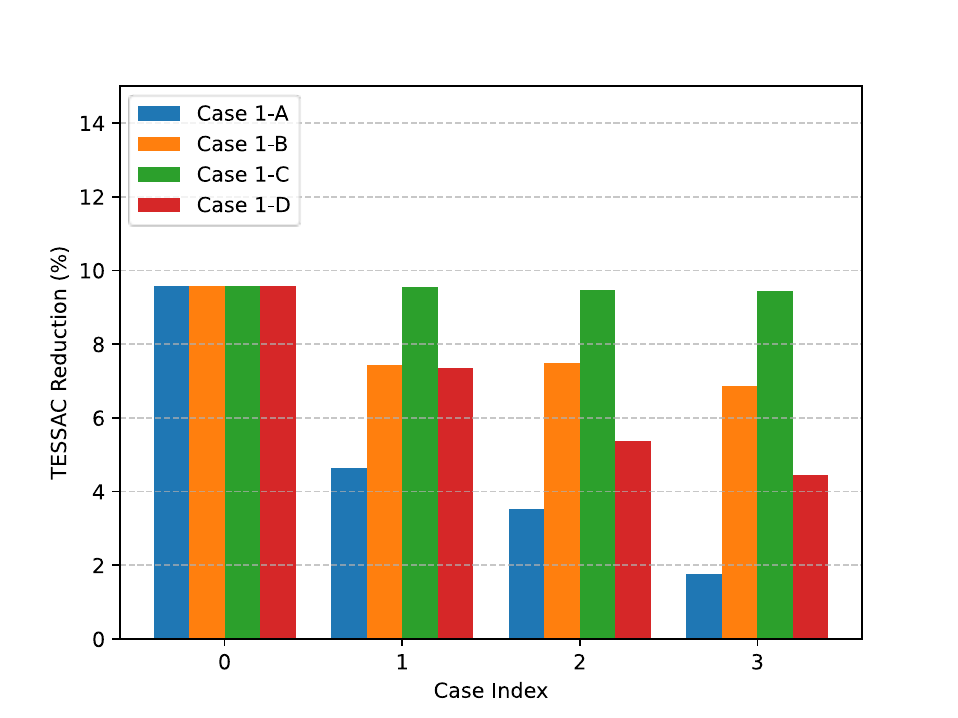}
\caption{TESSAC reduction obtained by joint replenishment strategy}
\label{fig: parametric study results}
\end{figure}

Figure \ref{fig: parametric study results} shows that differences in orbital altitude (Case C) have minimal impact on cost-effectiveness: a 300 km variation leads to only about a 0.1\% difference. By contrast, variations in other parameters---such as shipping size, the number of orbital planes, or the number of satellites per plane---cause a more noticeable reduction in the benefits offered by joint replenishment. Moreover, the degree of this reduction increases as the magnitude of the parameter differences.

\subsection{Case Study 2: Maintenance of Three Satellite Constellations}
The second case study focuses on three satellite constellations with parameters provided in Table \ref{tab: example constellation parameters}. From a maintenance perspective, these constellations share similar characteristics, including comparable constellation patterns, failure rates, and shipping sizes.
\begin{table}[hbt!]
    \centering
    \caption{Parameters for three satellite constellations for Case Study 2}
    \begin{tabular}{lccccc}
        \hline\hline
        Parameter & Notation & \multicolumn{3}{c}{Value} & Unit\\
        & & $j=1$ & $j=2$ & $j=3$ & \\\hline
        Altitude of orbital planes in constellation & $h_{\text{plane},j}$ & 1100 & 1300 & 1200 & km\\
        Satellite failure rate & $\lambda_{\text{sat},j}$ & 0.1 & 0.11 & 0.12 & satellites/year\\
        Number of orbital planes in constellation& $N_{\text{plane},j}$ & 24 & 26 & 20 & planes\\
        Number of satellites per orbital planes & $N_{\text{sat},j}$ & 20 & 22 & 24 & satellites/plane\\
        Shipping size of individual satellite & $v_{j}$ & 1 & 2 & 2 & \# slots\\
        Satellite dry mass & $M_{\text{dry},j}$ & 200 & 280 & 350 & kg\\
        Satellite manufacturing cost & $c_{\text{manufac},j}$ & 0.5 & 1 & 1 & \$M/satellite\\
        Mass flow rate of satellite's propulsion system & $\dot{M}_{\text{prop},j}$ & 1.7E-5 & 2.4E-5 & 3.0E-5 & kg/s\\
        Inclination & $i$ & \multicolumn{3}{c}{60} & deg\\
        Fuel mass conversion coefficient & $\epsilon_{\text{fuel},j}$ & \multicolumn{3}{c}{0.01} & \$M/kg\\
        Annual satellite holding cost & $c_{\text{hold},j}$ & \multicolumn{3}{c}{0.5} & \$M/satellite/year \\
        Exhaust velocity of satellite's propulsion system & $V_{\text{ex},j}$ & \multicolumn{3}{c}{11.77} & km/s\\
        Specific impulse of satellite's propulsion system & $I_{\text{sp},j}$ & \multicolumn{3}{c}{1200} & sec\\
        \hline\hline
    \end{tabular}
    \label{tab: example constellation parameters}
\end{table}

\begin{table}[hbt!]
    \centering
    \caption{Decision variables of independent replenishment strategy for Case Study 2}
    \begin{tabular}{lccc}
        \hline\hline
        Parameter & Notation & Range & Unit\\\hline
        ROP of in-plane spares & $s_1$ & $\left[1, 10\right]$ & satellites\\
        Order quantity of in-plane spares & $Q_1$ & $\left[1, Q_{\text{max}}\right]$ & satellites\\
        SROP of parking spares & $k_{R_1,3}$ & $\left[1, 40\right]$ & - \\
        Order quantity of parking spares & $k_{Q_1,3}$ & $\left[1, Q_{\text{max}}\right]$ & - \\
        Number of parking orbits & $N_{\text{parking}}$ & $\left[1, 20\right]$ & planes\\
        Altitude of parking orbits & $h_{\text{parking}}$ & $\left[500, 1000\right]$ & km\\
        \hline\hline
    \end{tabular}
    \label{tab: case study 2 independent decision variables}
\end{table}

Table \ref{tab: case study 2 independent decision variables} presents the decision variables for the independent replenishment strategy. Additionally, $\mathcal{H}_{\text{parking}}$ in Eq. \eqref{eqn: parking orbits values} is adopted to define the admissible altitude of the parking orbits. Table \ref{tab: reference solution} shows the optimal solutions of the independent replenishment strategy with each launch option.

The results show that, compared to the normal launcher, using the mega launcher generally reduces launch costs but increases holding costs. This trade-off yields different optimal solutions among the constellations. Specifically, constellations 2 and 3 benefit more from the mega launcher, whereas constellation 1 favors a normal launcher. Constellation 1 cannot fully exploit the lower satellite launch cost of the larger launcher because its high capacity necessitates holding additional spare satellites, negating the overall benefit of reduced launch costs.

Another notable observation is that using the mega launcher involves fewer parking orbits at lower altitudes. The reduced number of parking orbits, combined with the larger capacity of the mega launcher, makes each parking orbit hold more spare satellites. Moreover, the lower altitude shortens the alignment cycle between target and parking orbits, providing more frequent maneuver windows---albeit at the expense of increased fuel consumption to raise the orbit.

\begin{table}[hbt!]
    \centering
    \caption{Optimal solution and resultant parameters of the independent replenishment strategy for each launch option}
    \begin{tabular}{lccccccc}
        \hline\hline
         & \multicolumn{3}{c}{Normal Launcher} & & \multicolumn{3}{c}{Mega Launcher}\\
         & Con. 1 & Con. 2 & Con. 3 & & Con. 1 & Con. 2 & Con. 3\\\hline
        \textbf{Optimal Solution} &&&&&&&\\
        $s$, satellites & 2 & 2 & 2 & & 3 & 3 & 3 \\
        $Q$, satellites & 2 & 2 & 2 & & 5 & 3 & 3 \\
        $k_{s,\text{parking}}$, - & 5 & 6 & 4 & & 13 & 15 & 13 \\
        $k_{Q,\text{parking}}$, - & 20 & 10 & 10 & & 39 & 41 & 41 \\
        $N_{\text{parking}}$, planes & 3 & 4 & 5 & & 1 & 2 & 2 \\
        $h_{\text{parking}}$, km & 600 & 750 & 650 & & 500 & 750 & 600 \\\addlinespace
        \textbf{Result Parameter} & & \\
        $TESSAC$, \$M & 178.6 & 349.3 & 320.7 & & 191.7 & 297.8 & 268.9 \\
        $C_{\text{lau}}$, \$M & 80.4 & 210.8 & 193.0 & & 49.2 & 102.3 & 93.7 \\
        $C_{\text{hold}}$, \$M & 72.3 & 71.5 & 65.2 & & 115.9 & 128.5 & 112.3 \\
        $\rho_{\text{plane}}$, \% & 98.1 & 98.1 & 98.1 & & 98.4 & 98.5 & 98.1 \\
        $\rho_{\text{parking}}$, \% & 98.6 & 98.6 & 98.2 & & 98.0 & 98.3 & 98.1 \\
        \hline\hline
    \end{tabular}
    \label{tab: reference solution}
\end{table}

\subsubsection{Solution in Centralized Scenario}
In the centralized scenario, the optimal solution of $\mathbf{P}_{\text{JC}}$ is obtained and compared to the solution from the independent replenishment strategy. Only the mega launcher is considered, and the design space for the decision variables is provided in Table \ref{tab: joint strategy decision variables in centralized scenario}. Here, $\mathcal{H}_{\text{parking}}$ are defined in the same way as in the independent replenishment strategy, with the results presented in Table \ref{tab: optimal centralized joint solution}.

Table \ref{tab: optimal centralized joint solution} summarizes the optimal solution along with the resulting inventory parameters. The results are compared to the best independent strategies for the three constellations with their optimal launch options---normal launcher for constellation 1, and mega launcher for constellations 2 and 3---and the differences in costs are provided in parentheses. Under the best independent strategies, the total TESSAC across all constellations amounts to \$745.3M (=\$178.6M+\$297.8.4M+\$268.9M). In contrast, the joint strategy achieves a reduced cost of \$718.2M, representing an approximate 3.6\% decrease in annual maintenance costs. This demonstrates that the joint replenishment strategy outperforms the independent strategy in this instance.

Specifically, the total annual launch cost of the independent strategy is \$276.4M (=\$80.4M+\$102.3M+\$93.7M), which decreases to \$258.5M under the joint replenishment strategy. This seems to be mainly attributed to the fact that constellation 1 can use the mega launcher by sharing launch opportunities with the other constellations.

In terms of holding costs, constellations 2 and 3 experience reductions under the joint replenishment strategy, with constellation 2's holding cost decreasing from \$128.5M to \$96.4M, and constellation 3's from \$112.3M to \$107.5M. However, constellation 1 sees an increase in holding cost, rising from \$72.3M to \$96.1M. Despite this increase, it remains lower than the holding cost under the independent replenishment strategy with the mega launcher, which amounts to \$115.9M.

Overall, by sharing launch opportunities with the other constellations, Constellation 1 can leverage the mega launcher more efficiently, achieving significantly lower inventory costs compared to the independent replenishment strategy with the mega launcher, where high holding costs offset many of its benefits. Additionally, the other constellations also benefit from reduced overall inventory costs by efficiently utilizing the mega launcher's large capacity through shared launch opportunities, even though they could employ the mega launcher independently as well.

\begin{table}[hbt!]
    \centering
    \caption{Decision variables of the joint replenishment strategy in the centralized scenario for Case Study 2, $j=1,2,3$}
    \begin{tabular}{lccc}
        \hline\hline
        Variables & Notation & Range & Unit\\\hline
        ROP of in-plane spares & $s_j$ & $\left[1, 10\right]$ & satellites\\
        Order quantity of in-plane spares & $Q_j$ & $\left[1, Q_{\text{max},j}\right]$ & satellites\\
        OUL of parking spare & $S_j$ & $\left[1, 40\right]$ & - \\
        SROP of parking spare & $U$ & $\left[200, 250\right]$ & \# slots\\
        The number of parking orbits & $N_{\text{parking}}$ & $\left[1, 20\right]$ & planes\\
        Altitude of parking orbits & $h_{\text{parking}}$ & $\left[500, 1000\right]$ & km\\
        \hline\hline
    \end{tabular}
    \label{tab: joint strategy decision variables in centralized scenario}
\end{table}

\begin{table}[hbt!]
    \centering
    \caption{Optimal solution and the resultant parameters of the joint replenishment strategy in the centralized scenario for Case Study 2}
    \begin{tabular}{lccc}
        \hline\hline
         & Con. 1 & Con. 2 & Con. 3\\\hline
        \textbf{Optimal Solution} &&&\\
        $s_j$, satellites & 3 & 3 & 3  \\
        $Q_j$, satellites & 5 & 5 & 10  \\
        $S_j$, - & 30 & 32 & 16 \\
        $U$, \# slots & \multicolumn{3}{c}{244}  \\
        $N_{\text{parking}}$, planes & \multicolumn{3}{c}{1}  \\
        $h_{\text{parking}}$, km & \multicolumn{3}{c}{500} \\\addlinespace
        \textbf{Result Parameter} & & \\
        $TESSAC$, \$M & \multicolumn{3}{c}{718.2 ($\downarrow 27.1$)} \\
        $C_{\text{lau}}$, \$M & \multicolumn{3}{c}{258.5 ($\downarrow 17.9$)} \\
        $C_{\text{hold},j}$, \$M & 96.1 ($\uparrow 23.8$) & 96.4 ($\downarrow 32.1$) & 107.5 ($\downarrow 4.8$) \\
        $\rho_{\text{plane},j}$, \% & 98.4 & 98.2 & 98.2 \\
        $\rho_{\text{parking},j}$, \% & 98.1 & 98.3 & 98.3\\
        \hline\hline
    \end{tabular}
    \label{tab: optimal centralized joint solution}
\end{table}

\subsubsection{Solution in Decentralized Scenario}
In the decentralized scenario, we solve $\mathbf{P}_{\text{JD}}^1$ first. The decision variables are listed in Table \ref{tab: joint strategy decision variables in decentralized scenario}. The admissible set of values for the altitude of the parking orbits, $\mathcal{H}_{\text{parking}}$, is defined in the same way as in the independent strategy. Additionally, we set $TESSAC_{\text{ref},j}$ to the TESSAC of the best independent strategy, with values of \$178.6M, \$297.8M, and \$268.9M for $j=1$, $2$, and $3$, respectively.
\begin{table}[hbt!]
    \centering
    \caption{Decision variables of the joint replenishment strategy in the decentralized scenario for Case Study 2, $j=1,2,3$}
    \begin{tabular}{lccc}
        \hline\hline
        Variables & Notation & Range & Unit\\\hline
        ROP of in-plane spares & $s_j$ & $\left[1, 10\right]$ & satellites\\
        Order quantity of in-plane spares & $Q_j$ & $\left[1, Q_{\text{max},j}\right]$ & satellites\\
        OUL of parking spares & $S_j$ & $\left[1, 40\right]$ & - \\
        SROP of parking spares & $U$ & $\left[200, 250\right]$ & \# slots\\
        The number of parking orbits & $N_{\text{parking}}$ & $\left[1, 20\right]$ & planes\\
        Altitude of parking orbits & $h_{\text{parking}}$ & $\left[500, 1000\right]$ & km\\
        Launch cost payment fraction & $y_j$ & $\left[0, 1\right]$ & -\\
        \hline\hline
    \end{tabular}
    \label{tab: joint strategy decision variables in decentralized scenario}
\end{table}

\begin{figure}[hbt!]
    \centering
    \includegraphics[page=6, width=1.0\textwidth]{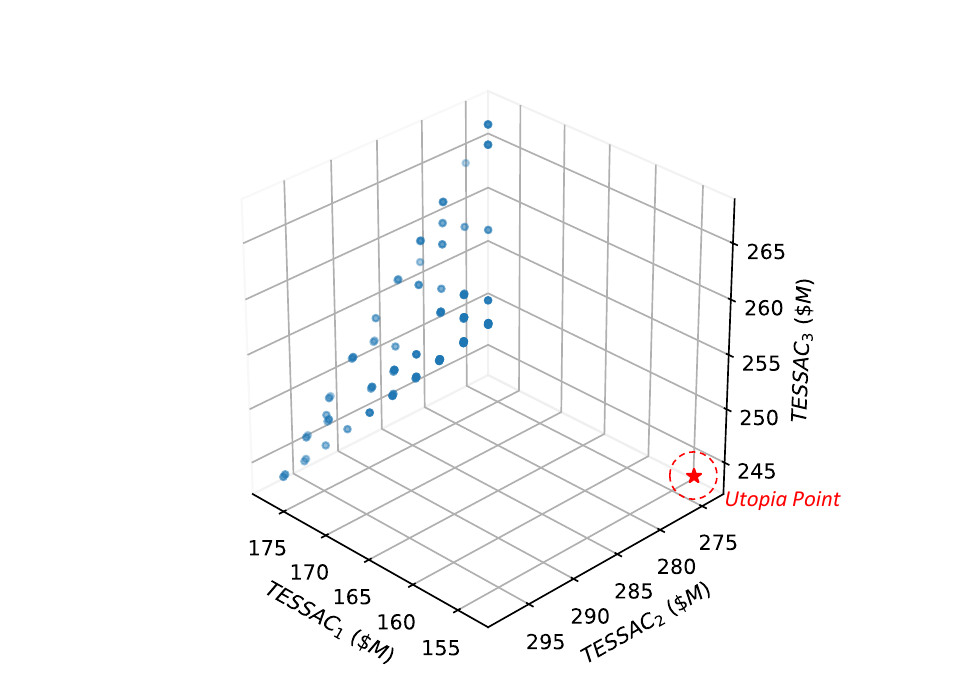}
    \caption{Efficient solutions of $\mathbf{P}_{\text{JD}}^1$ in the objective space in the decentralized scenario for Case Study 2}
\label{fig: Pareto solutions}
\end{figure}

With this decision space and solver, we obtain 100 solutions. Figure \ref{fig: Pareto solutions} presents the resulting solutions in the objective space. The red star in the figure represents the \textit{utopia point}, which lies in the decreasing direction of all three objectives. Although it is unattainable, it serves as a guide towards the desired optimization direction.

Next, based on the obtained set of efficient solutions, we solve $\mathbf{P}_{\text{JD}}^2$. The parameters reflecting each constellation operator’s influence on the final solution are then set accordingly. Specifically, we consider two different parameter sets $\left(\beta_1,\beta_2,\beta_3\right)$: $\left(0.2,0.4,0.4\right)$ for the first instance and $\left(0.6,0.2,0.2\right)$ for the second instance.

Table \ref{tab: selected decentralized joint solution} presents the solutions obtained by solving $\mathbf{P}_{\text{JD}}^2$ for each instance. In the first instance, all decision variables, except for the launch cost payment fraction, are the same as those in the centralized scenario. This means that the maintenance system can fully leverage the benefits of the joint replenishment strategy with this solution. With launch cost payment fractions of 0.21, 0.47, and 0.32 for constellations 1, 2, and 3, respectively, constellations 1 and 3 pay less for launch services than in the independent strategy, while constellation 2 pays more. However, the higher payment by constellation 2 is offset by reduced holding costs. As a result, all constellations benefit from the joint replenishment strategy with the selected solution.

In the second instance, SROP is slightly different from the solution of the centralized scenario, but the other decision variables remain the same. Although the replenishment strategy is similar, compared to the first instance, this instance places more weight on the costs for constellation 1's operator. As a result, constellation 1's TESSAC is reduced, while the TESSACs for the other two operators increase. This is mainly attributed to the launch cost payment fractions, which are 0.12, 0.51, and 0.37 for constellations 1, 2, and 3, respectively. While the value for constellation 1 has increased compared to the first instance, the values for constellations 2 and 3 have decreased compared to the first instance.

\begin{table}[hbt!]
    \centering
    \caption{Selected solutions and the resultant parameters of the joint replenishment strategy in the decentralized scenario for Case Study 2}
    \begin{tabular}{lccccccc}
        \hline\hline
        & \multicolumn{3}{c}{First Inst.} && \multicolumn{3}{c}{Second Inst.} \\
        & Con. 1 & Con. 2 & Con. 3 && Con. 1 & Con. 2 & Con. 3\\\hline
        \textbf{Bargaining Power} &&&&&&&\\
        $\beta_j$, -& 0.2 & 0.4 & 0.4 && 0.6 & 0.2 & 0.2 \\
        \textbf{Optimal Solution} &&&&&&&\\
        $s_j$, satellites & 3 & 3 & 3 && 3 & 3 & 3  \\
        $Q_j$, satellites & 5 & 5 & 10 && 5 & 5 & 10 \\
        $S_j$, - & 30 & 32 & 16 && 30 & 32 & 16 \\
        $U$, \# slots & \multicolumn{3}{c}{244} && \multicolumn{3}{c}{247} \\
        $N_{\text{parking}}$, planes & \multicolumn{3}{c}{1} && \multicolumn{3}{c}{1} \\
        $h_{\text{parking}}$, km & \multicolumn{3}{c}{500} && \multicolumn{3}{c}{500} \\
        $y_j$ (-) & 0.21 & 0.47 & 0.32 && 0.12 & 0.51 & 0.37 \\
        \addlinespace
        \textbf{Result Parameter} &&&&&&&\\
        $TESSAC_j$, \$M & 177.0 & 287.1 & 254.1 && 153.4 & 297.7 & 267.2 \\
        $y_jC_{\text{lau}}$, \$M & 54.3 & 121.5 & 82.7 && 31.2 & 132.8 & 96.3\\
        $C_{\text{hold},j}$, \$M & 96.1 & 96.4 & 107.5 && 95.6 & 95.9 & 107.0 \\
        $\rho_{\text{plane},j}$, \% & 98.4 & 98.2 & 98.2 && 98.4 & 98.3 & 98.2 \\
        $\rho_{\text{parking},j}$, \% & 98.1 & 98.3 & 98.3 && 98.0 & 98.2 & 98.2 \\
        \hline\hline
    \end{tabular}
    \label{tab: selected decentralized joint solution}
\end{table}

\section{Conclusions}\label{sec: section 6}
This work proposed a novel joint replenishment strategy for multiple satellite constellations based on the supply chain of spare satellites with shared launch opportunities and parking orbits. Following this supply chain, an inventory management model with a parametric policy---$\left(U,\mathbf{S}\right)$ for parking spares and $\left(s,Q\right)$ for in-plane spares---was developed. 

Using this model, we proposed solution frameworks for two different scenarios---centralized and decentralized. In the centralized scenario, where a single authority determines all parameters for the joint replenishment strategy, a single-objective optimization problem is solved to minimize the total annual expenditure of the entire maintenance system. In the decentralized scenario, multiple steps are required to reach an agreement.

We conducted two case studies to examine how parameter dissimilarities among constellations affect the cost-effectiveness of the joint replenishment strategy---generally reducing its overall benefits---and to demonstrate the proposed solution procedure in a practical context.

This work can be extended in various ways. One promising avenue is to explore the interactions among different operators in greater depth. For instance, an operator might deviate from an agreement to pursue its own objectives. To prevent such deviations and maintain maximum efficiency in the overall supply chain, it is crucial to investigate appropriate contract structures in this cooperative scenario, particularly considering the unique aspects of maintaining satellite constellations.

Another interesting direction is to investigate scenarios where substitution between constellations is permitted. For instance, if two constellations operate satellites with similar missions and payloads but differ in capacity, smaller-capacity satellites could potentially substitute for a larger-capacity one, or vice versa. This substitution mechanism could introduce new dynamics and optimization opportunities in the joint replenishment strategy.

\section*{Acknowledgments}
This work was supported by the National Research Foundation of Korea (NRF) grant funded by the Korea government (MSIT) (RS-2025-00563732). The authors thank the NRF for the support of this work.

\end{document}